\documentclass[11pt,a4paper]{amsart}

\usepackage{amssymb}
\usepackage[all]{xypic}
\swapnumbers
\def\frak{\mathfrak}
\def\Bbb{\mathbb}
\def\Cal{\mathcal}

\newtheorem{prop}[subsection]{Proposition}
\newtheorem*{prop*}{Proposition}
\newtheorem{thm}[subsection]{Theorem}
\newtheorem*{thm*}{Theorem}
\newtheorem{lem}[subsection]{Lemma}
\newtheorem*{lem*}{Lemma}
\newtheorem{kor}[subsection]{Corollary}
\newtheorem*{kor*}{Corollary}
\newtheorem{rem}[subsection]{Remark}

\newcommand{\ad}{\operatorname{ad}}
\newcommand{\Ad}{\operatorname{Ad}}
\renewcommand{\exp}{\operatorname{exp}}
\newcommand{\id}{\operatorname{id}}
\newcommand{\Ker}{\operatorname{Ker}}
\renewcommand{\ker}{\operatorname{ker}}
\newcommand{\im}{\operatorname{im}}
\newcommand{\Hom}{\operatorname{Hom}}

\newcommand{\Alt}{\operatorname{Alt}}

\newcommand{\x}{\times}
\renewcommand{\o}{\circ}
\newcommand{\ddt}{\tfrac{d}{dt}|_{t=0}}

\let\ccdot\cdot
\def\cdot{\hbox to 2.5pt{\hss$\ccdot$\hss}}

\newcommand{\al}{\alpha}

\newcommand{\ga}{\gamma}
\newcommand{\io}{\iota}
\newcommand{\ka}{\kappa}
\newcommand{\la}{\lambda}
\newcommand{\om}{\omega}
\renewcommand{\phi}{\varphi}
\newcommand{\ph}{\varphi}
\newcommand{\ps}{\psi}
\newcommand{\si}{\sigma}

\newcommand{\ze}{\zeta}
\newcommand{\De}{\Delta}
\newcommand{\Ga}{\Gamma}
\newcommand{\La}{\Lambda}
\newcommand{\Om}{\Omega}
\newcommand{\Ph}{\Phi}
\newcommand{\Ps}{\Psi}
\newcommand{\Si}{\Sigma}

\def\sideremark#1{\ifvmode\leavevmode\fi\vadjust{
\vbox to0pt{\hbox to 0pt{\hskip\hsize\hskip1em
\vbox{\hsize3cm\tiny\raggedright\pretolerance10000
\noindent #1\hfill}\hss}\vbox to8pt{\vfil}\vss}}}

\begin{document}
\title{Bernstein--Gelfand--Gelfand
sequences} 
\author{Andreas \v Cap, Jan Slov\'ak, and Vladim\'{i}r Sou\v cek}
\address{A.C.: Institut f\"ur Mathematik, Universit\"at Wien, Strudlhofgasse 4,
A--1090 Wien, Austria and International Erwin Schr\"odinger Institute for
Mathematical Physics, Boltzmanngasse 9, A--1090 Wien, Austria\newline\indent
J.S.: Department of Algebra and Geometry, Masaryk University in Brno,
Jan\'a\v ckovo n\' am.\ 2a, 662~95 Brno, Czech
Republic\newline\indent
V.S.: Mathematical Institute, Charles University, Sokolovsk\'a 83,
Praha, Czech Republic}
\email{Andreas.Cap@esi.ac.at, slovak@math.muni.cz,\newline\indent soucek@karlin.mff.cuni.cz}
\begin{abstract}
This paper is devoted to the study of geometric structures modeled on
homogeneous spaces $G/P$, where $G$ is a real or complex semisimple
Lie group and $P\subset G$ is a parabolic subgroup. We use methods from
differential geometry and very elementary finite--dimensional
representation theory to construct sequences of invariant differential
operators for such geometries, both in the smooth and the holomorphic
category. For $G$ simple, these sequences specialize on the homogeneous 
model $G/P$ to the celebrated (generalized)
Bernstein--Gelfand--Gelfand resolutions in the holomorphic category,
while in the smooth category we get smooth analogs of these
resolutions. In the case of geometries locally isomorphic to the
homogeneous model, we still get resolutions, whose cohomology is
explicitly related to a twisted de~Rham cohomology. In the general
(curved) case we get distinguished curved analogs of all the invariant
differential operators occurring in Bernstein--Gelfand--Gelfand
resolutions (and their smooth analogs). 

On the way to these results, a significant part of the general theory
of geometrical structures of the type described above is presented
here for the first time. 
\end{abstract}

\maketitle

\section{Introduction}\label{1}

Our approach to geometries modeled on homogeneous spaces goes back to
E.~Cartan's notion of an `espace generalis\'e'. The central
objects for such geometries are suitably normalized Cartan connections
in the sense commonly adopted, see e.g. \cite{Sharpe}. The models for
the geometries considered in this paper are homogeneous spaces of the
type $G/P$, where $G$ is real or complex semisimple and $P\subset G$
is a parabolic subgroup. In this case, there is a close link to the
project of parabolic invariant theory suggested by Ch.~Fefferman in
\cite{Fef} and in view of this context we talk about the (real and
complex) {\em parabolic geometries}. 

We explore the semi--holonomic jet modules and we use implicitly the
cohomological information given by Kostant's version of the
Bott--Borel--Weil theorem in order to construct sequences of
homomorphisms between jet--modules, which in turn give rise to
sequences of invariant differential operators expressed in terms of
the invariant derivatives with respect to Cartan connections, on all
(curved) geometries in question.  These sequences are differential
complexes if certain twisted de~Rham sequences are complexes, and then
they compute the same cohomology. In particular, this always happens
for the homogeneous models themselves and then our sequences
specialize to the Bernstein--Gelfand--Gelfand resolutions well known
from representation theory for complex $G/P$, while their real smooth
analogues are provided for all real forms of this situation.

In spite of the fact that we have mentioned a few concepts from
representation theory, we want to underline that no deeper
aspects of representation theory are used in the construction of our new
sequences of invariant operators and in the discussion of their basic
properties. In particular, no infinite dimensional representation theory
is needed. It is rather the language and the way of thinking of
representation theory that is essential (in a similar way as the
categorical language is useful in mathematics even if no deep results
of category theory are used). In order to stress this feature, we have
postponed the more detailed analysis of the structure of the sequences
to a forthcoming second part of the article and we hope that the first
part is accessible for differential geometers without a deeper
background in representation theory. We also provide a quite detailed 
exposition of the necessary algebraic background. In particular
we have included two appendices covering some material which is rather
well known in representation theory.

The first general geometric theory close to our needs
had been worked out in the series of papers by N.~Tanaka and his school
aiming at the original equivalence problem of E.~Cartan, 
see \cite{Tan, Yam, Mor} and the
references therein. Our inspiration comes, however, rather from the interest
in the links between twistor theory and representation theory, as
explained in the book \cite{BE}. In the generality we need, the normalized
Cartan connections were constructed in \cite{CS} first. We have been also
influenced by the translation principle in representation theory (see
\cite{BoeColl1, BoeColl2} for example) and, in particular, by some ideas in
the second part of Baston's paper \cite{Bas}. Some arguments and
proofs in the latter paper seem very unclear to us, however.

There are also many treatments of specific examples of parabolic
geometries in the literature, including e.g.\ projective, conformal,
almost Grassmannian, and CR--geometries. Most of these well known geometries
correspond to the so called $|1|$--graded Lie algebras $\frak g$ which can be
equivalently expressed by the requirement that the tangent spaces
correspond to irreducible representations of the parabolic subgroup
$P$. Our theory of semi--holonomic jet--modules is in fact a
generalization of the approach worked out for all real $|1|$--graded
algebras in our former papers \cite{CSS1, CSS2, CSS3} (and this paper could
be also viewed as a fourth part of this series expanded to the full
generality of parabolic geometries).  On the other hand, there are only few
explicit examples of curved analogues of the Bernstein--Gelfand--Gelfand
resolutions available in the literature, see e.g. \cite{MikeBGG}, and in
fact only the case of conformal Riemannian geometries has been studied
systematically, see \cite{Gov} and \cite{ES} for two different approaches. 
For an introduction addressed to wide audience, see the forthcoming paper
\cite{MikeNotices}.

Let us indicate the structure of the paper. In the next section, we first
collect the necessary information on $|k|$-graded Lie algebras and the
structure of the corresponding Lie groups, and then real and complex
parabolic geometries are introduced (cf. \ref{2.8}).  Our point of view is
that the geometry on a manifold $M$ is given by a \emph{choice} of a
Cartan connection (with possible further normalization) and we are
interested in the general calculus which such a choice offers. In a certain
sense, this is similar to the r\^ole of the general calculus for
linear connections in Riemannian geometry by application to the
Levi--Civita connection. Thus we only briefly discuss the more
classical underlying geometrical information on the manifolds $M$
themselves and the question of constructing a (normalized) Cartan
connection from these more basic data, cf. \ref{2.13}.  See \cite{CS,
Mor} for more information on this aspect. We also introduce the
concepts of natural bundles and operators for parabolic geometries in
the end of Section \ref{2}. 

The third section deals with our basic algebraic tool, the semi--holonomic
jet modules. The \emph{invariant derivative} with respect to Cartan
connections then leads to the notion of \emph{strongly invariant
differential operators}
which are defined by means of $P$--module homomorphisms. As a first
application, we introduce the twisted exterior derivatives which are certain
torsion adjusted versions of the covariant exterior derivatives
induced by the Cartan connections on certain bundles.

The main results are stated and proved in Section \ref{4}.  Referring
implicitly to the structure of the Lie algebra cohomologies, we first embed
the natural vector bundles corresponding to cohomologies into exterior forms
by means of distinguished differential operators $L$, see Theorem \ref{5.6}. 
Then we use the twisted exterior derivatives in order to construct
explicitly many $P$--module homomorphisms of the semi--holonomic jet
modules, cf. Proposition \ref{5.5a}. The corresponding invariant
differential operators build the
\emph{Bernstein--Gelfand--Gelfand sequences}. Finally we discuss the conditions
under which these sequences form differential complexes, and we discuss their
cohomologies, cf. \ref{5.8}--\ref{5.9a}.

Finally, we illustrate briefly the achievements
on at least one non--trivial parabolic geometry and this is done in
Section \ref{5}.

Throughout the paper, we discuss the real and complex manifolds and groups
at the same time. We should point out however, that the relation 
between the real and complex settings deserves more attention. In fact, we
are able to present both smooth and holomorphic results in one line of
arguments, because our point is to use the $P$--module homomorphisms in order
to construct the sequences of operators. The distinction is hidden in the
explicit structure of the Lie algebra cohomologies, which we use only
implicitly. One should say, however, this {\it does not mean\/} that working
out the details for one real form gives explicit results for all other real
or complex forms of the group in question. This ambiguity disappears
only if we restrict ourselves to complex representations of the real forms.

A more detailed discussion of our Bernstein--Gelfand--Gelfand sequences
requires a deeper study of the cohomological information. Essentially, the
non--trivial operators between the irreducible bundles in the sequence
correspond to arrows in the Hasse diagram of the parabolic subalgebras and
the knowledge of this structure leads to quite explicit information on the
individual operators. We have preferred to postpone all considerations
which need more involved information from representation theory to a
prospective continuation in order to keep the flavor of this article.  

\noindent{\bf Acknowledgements.} The research evolved during a 
stay of the first two authors at the University of Adelaide supported by the
Australian Research Council, and during the meetings of all three authors at
the Erwin Schr\"odinger Institute for Mathematical Physics in Vienna, the
Masaryk University in Brno, and the Charles University in Prague. The
institutional support by GA\v CR, Grant Nr.~201/99/0675 has been essential,
too. Our particular thanks are due to Michael Eastwood who explained
to us several aspects of the Bernstein--Gelfand--Gelfand resolutions.

\section{Parabolic geometries}\label{2}
In this section we review basic facts about $|k|$--graded Lie algebras
and we give basic definitions on parabolic geometries and invariant
differential operators on manifolds equipped with geometries of that
type. Most of the facts on the algebras go back to \cite{Tan, Yam}, see also 
\cite{CS} which is fully compatible in notation. 

\subsection{Definition}\label{2.1} Let $\Bbb K$ be $\Bbb R$ or $\Bbb C$. 
A $|k|$--graded Lie algebra over $\Bbb K$, $k\in \Bbb N$ 
is a Lie algebra $\frak g$ over $\Bbb K$ together with a decomposition 
$$
\frak g=\frak g_{-k}\oplus\dots
\oplus \frak g_{-1}\oplus\frak g_0\oplus\frak g_1\oplus\dots\oplus \frak g_k
$$
such that $[\frak g_i,\frak g_j]\subset\frak g_{i+j}$ and such that the 
subalgebra $\frak g_-:=\frak g_{-k}\oplus\dots\oplus\frak g_{-1}$ is 
generated by $\frak g_{-1}$. In the whole paper, we will only
deal with {\em semisimple\/} $|k|$--graded Lie algebras.

By $\frak p$ we will denote the subalgebra 
$\frak g_0\oplus\dots\oplus\frak g_k$
of $\frak g$, and by $\frak p_+$ the subalgebra 
$\frak g_1\oplus\dots\oplus\frak g_k$ of $\frak p$.

There is always a unique element $E\in\frak g$ whose adjoint action
is given by $[E,X]=\ell X$ for $X\in\frak g_\ell$. The element $E$ is
contained in the center of the subalgebra $\frak g_0$, which is
always reductive. Using this, one shows that any ideal of $\frak g$ is
homogeneous. Thus, a semisimple $|k|$--graded Lie algebra is always a
direct sum of simple $|k_i|$--graded Lie algebras, where all $k_i\leq
k$. Hence, one usually can reduce most discussions to the simple
case. When dealing with the semisimple case, we have to assume that
none of the simple factors is contained in $\frak g_0$, for technical 
reasons. Since
basically we are interested in homogeneous spaces $G/P$, where $G$ is
a Lie group with Lie algebra $\frak g$ and $P$ an appropriate subgroup
with Lie algebra $\frak p$, and their curved analogs, this is not
really a restriction. 

For each $i=1,\dots, k$, the Killing form
of $\frak g$ induces an isomorphism $\frak g_{i}\cong\frak g_{-i}^*$
of $\frak g_0$--modules. Finally, the powers of $\frak p_+$ are given
by $\frak p_+^i=\frak g_i\oplus\dots\oplus\frak g_k$, for
$i=1,\dots,k$. See e.g. \cite[Section 3]{Yam} for details. 

\subsection{}\label{2.2}
In the complex case, the meaning of a $|k|$--grading is particularly simple
to describe. One can show that there always exists a Cartan subalgebra
$\frak h\subset\frak g$ which contains the element $E$ from above, and a
choice of positive roots $\De_+$ for $\frak h$ such that all root spaces
corresponding to simple roots are either contained in $\frak g_0$ or in
$\frak g_1$. Denoting by $\Si$ the set of those simple roots, whose root
spaces are contained in $\frak g_1$, one sees that the grading on $\frak g$
is given by the $\Si$--height of roots. That is, if $\al$ is a root, then
the root space $\frak g_\al$ is contained in $\frak g_i$, where $i$ is the
sum of all coefficients of elements of $\Si$ in the expansion of $\al$ as a
linear combination of simple roots. In particular, this implies that the
subalgebra $\frak p$ is always a parabolic subalgebra of $\frak g$, and
$\frak p=\frak g_0\oplus\frak p_+$ is exactly the Levi decomposition of
$\frak p$ into a reductive and a nilpotent part.  

Conversely, if $\frak g$ is complex and semisimple and $\frak p\subset\frak
g$ is a parabolic subalgebra, then one can find a Cartan subalgebra and a
set of positive roots such that $\frak p$ is the standard parabolic
corresponding to a set $\Si$ of simple roots. But then the $\Si$--height as
defined above gives a $|k|$--grading on $\frak g$, where $k$ is the
$\Si$--height of the maximal root of $\frak g$, such
that $\frak p=\frak g_0\oplus\dots\oplus\frak g_k$.  See e.g. \cite[p.
88]{Humphreys} or \cite[Section 2]{BE} for more details.

Thus, in the complex case giving a $|k|$--grading on $\frak g$ is the
same thing as giving a parabolic subalgebra $\frak p$ of $\frak
g$. Therefore, complex $|k|$--graded semisimple Lie algebras can be
conveniently denoted by Dynkin diagrams with crossed nodes. That is,
given a $|k|$--graded semisimple complex Lie algebra we may assume
that $\frak p$ is the standard parabolic subalgebra corresponding to a
set $\Si$ of simple roots. Then we denote the $|k|$--graded Lie
algebra $\frak g$ by crossing out the nodes corresponding to the
simple roots contained in $\Si$ in the Dynkin diagram of $\frak
g$. See the book \cite{BE} for a detailed discussion of the Dynkin
diagram notation for parabolic subalgebras. 

Finally note that for a $|k|$--graded Lie algebra $\frak g$ over $\Bbb
R$ the complexification $\frak g^\Bbb C$ of $\frak g$ is
$|k|$--graded, too. So in general we deal with certain real forms of
pairs $(\frak g,\frak p)$, where $\frak g$ is complex and semisimple
and $\frak p$ is a parabolic in $\frak g$. The classification of all these
real forms is provided in \cite[Section 4]{Yam}.

\subsection{}\label{2.3}
Suppose that $\frak g$ is $|k|$--graded and semisimple over $\Bbb
K=\Bbb R$ or $\Bbb C$, and let $G$ be any Lie group with Lie algebra
$\frak g$. (We do not assume that $G$ is connected.) Then we can
define subgroups $G_0\subset P\subset G$ as follows: $G_0$ consists of
all elements of $G$ such that the adjoint action $\Ad(g):\frak
g\to\frak g$ of $g$ preserves the grading of $\frak g$. By $P$ we
denote the subgroup of all elements $g\in G$ such that $\Ad(g)$
preserves the filtration by right ends induced by the grading of
$\frak g$, i.e.\ {}$\Ad(g)(\frak g_i)\subset\frak
g_i\oplus\dots\oplus\frak g_k$. By definition $G_0$ is a subgroup of
$P$, and one easily verifies that $G_0$ and $P$ have Lie algebras
$\frak g_0$ and $\frak p$, respectively, see
e.g. \cite[2.9]{CS}. Moreover, it can be shown that if $\frak g$ is
simple, then $P$ equals the normalizer $N_G(\frak p)$ of $\frak p$ in
$G$, so it is the usual parabolic subgroup associated to the parabolic
subalgebra $\frak p$. 

The following proposition clarifies the structure of the group $P$: 
\begin{prop*}
Let $g\in P$ be any element. Then there exist unique elements 
$g_0\in G_0$ and $X_i\in\frak g_i$ for $i=1,\dots, k$, such that 
$$g=g_0\exp(X_1)\dots\exp(X_k).$$
\end{prop*}
\begin{proof}
See \cite[2.10]{CS}.
\end{proof}

\subsection{}\label{2.4}
For $i=1,\dots,k$ we define a subgroup $P_+^i\subset P$ as the image
under the exponential map of $\frak g_i\oplus\dots\oplus\frak g_k$,
and we write $P_+$ for $P_+^1$. Then we have $P\supset P_+\supset
P_+^2\supset\dots\supset P_+^k$. The subgroup $P_+\subset P$ is
obviously normal and by Proposition \ref{2.3} we have $P/P_+\cong
G_0$, so $P$ is the semidirect product of $G_0$ and the normal nilpotent
subgroup $P_+$. More generally, for each $i>1$ we see that $P/P_+^i$
is the semidirect product of $G_0$ and the normal nilpotent subgroup
$P_+/P_+^i$. 

The adjoint action of $P$ on $\frak g$ by definition preserves any of
the subspace $\frak g_i\oplus\dots\oplus\frak g_k$ for $i=-k,\dots,k$.
Thus for each $i=-k,\dots, k$ and $j>i$ we get an induced action of
$P$ on the quotient $\frak g_i\oplus\dots\oplus\frak g_k/(\frak
g_j\oplus\dots\oplus\frak g_k)$. With a slight abuse of notation, we will
denote this $P$--module by $\frak g_i\oplus\dots\oplus\frak g_{j-1}$. Again
by Proposition \ref{2.3}, the action of $P_+^{j-i}$ on
$\frak g_i\oplus\dots\oplus\frak g_{j-1}$ is trivial, so the action of
$P$ on this space is induced by an action of $P/P_+^{j-i}$. In
particular, we get an action of $P$ on $\frak g_-=\frak g/\frak p$,
which is induced by an action of $P/P_+^k$. 

There is another important consequence of Proposition \ref{2.3}:
Suppose that $\Bbb V$ and $\Bbb W$ are $P$--modules and that $\Ph:\Bbb V\to
\Bbb W$ is a
linear mapping. Suppose that $\Ph$ is equivariant for the action of
$G_0$ and for the (infinitesimal) action of $\frak g_1$. Since $\frak
p_+$ is generated by $\frak g_1$ this implies equivariancy with
respect to $\frak p_+$ and thus also with respect to $P_+$, so using
Proposition \ref{2.3} we see that $\Ph$ is actually a homomorphism of
$P$--modules. This will be technically very important in the sequel. 

\subsection{}\label{2.5}
For a Lie group $G$ with $|k|$--graded semisimple Lie algebra $\frak
g$ and the subgroup $P$ defined in \ref{2.3} above, consider the
homogeneous space $G/P$. This homogeneous space is the flat model for
the parabolic geometry of the type $(G,P)$ that we are going to
study. It is well known that the canonical projection $G\to G/P$ is a
principal fiber bundle with structure group $P$. 

If $G$ is a complex Lie group, then $P$ is a parabolic
subgroup, so $G/P$ is a generalized flag manifold, and thus in
particular a compact complex manifold. In the real case, $G/P$ need
not be compact in general, as the example of the conformal spheres in
indefinite signature shows.

Next suppose that $\la:P\to \text{GL}(\Bbb V)$ is a representation of $P$
on a finite dimensional vector space $\Bbb V$. Then we can form the
associated bundle $V:=G\x_P\Bbb V\to G/P$. This is a \emph{homogeneous
vector bundle}, that is the canonical left action of $G$ on $G/P$ lifts
to a left action of $G$ on $V$ by vector bundle
homomorphisms. Conversely, given a homogeneous vector bundle $E\to
G/P$, consider the fiber $\Bbb E$ of $E$ over the canonical base point
$o\in G/P$. Since the action of any element of $P$ on $G/P$ maps $o$
to itself, the action of $G$ on $E$ induces a representation of $P$ on
$\Bbb E$ and one easily verifies that $G\x_P\Bbb E$ and $E$ are
isomorphic homogeneous vector bundles (i.e.\ there is a
$G$--equivariant isomorphism of vector bundles between
them). Consequently, there is a bijective correspondence between
finite dimensional representations of $P$ and homogeneous vector
bundles over $G/P$. In the case where $G$ is a complex Lie group, the
action of $G$ on $G/P$ is holomorphic and there is a bijective
correspondence between holomorphic finite dimensional representations
of $P$ and holomorphic homogeneous vector bundles over $G/P$ (that is
holomorphic bundles with holomorphic $G$--actions). 

In particular, the tangent and cotangent bundles of $G/P$ are
homogeneous vector bundles. One easily verifies that they correspond
to the representations of $P$ on $\frak g_-\cong\frak g/\frak p$ and
$\frak p_+$ induced by the adjoint action, respectively. 
In the complex case, these
representations induce the holomorphic tangent and cotangent bundle. 

For a homogeneous vector bundle $E\to G/P$ consider the space $\Ga(E)$
of smooth sections of $E$. There is an induced action of $G$ on this
space given by $(g\cdot s)(x)=g\cdot (s(g^{-1}\cdot x))$ for $x\in
G/P$. In the complex case, we can deal similarly with the spaces of
holomorphic sections. 

\subsection*{Definition} Let $E$ and $F$ be homogeneous vector bundles
over $G/P$. A (linear) {\em invariant differential operator\/} 
$D:\Ga(E)\to\Ga(F)$ is a linear differential operator $D$ which is
equivariant for the $G$--actions constructed above. 

\subsection{}\label{2.6}
If $D$ is
of order $\leq r$, then it is induced by a vector bundle homomorphism
$\tilde D:J^r(E)\to F$, where $J^r(E)$ is the $r$--th jet prolongation of
$E$. Now simply by functoriality of the $r$--th jet prolongation,
$J^r(E)$ is again a homogeneous vector bundle, and the invariance of
$D$ is equivalent to the fact that $\tilde D$ is equivariant for the
$G$--actions on $J^r(E)$ and $F$. Since $G$ acts transitively on
$G/P$, the homomorphism $\tilde D$ is actually determined by its
restriction $\tilde D:J^r(E)_o\to F_o$ to the fiber over $o\in G/P$,
and by invariance of $D$, this map is $P$--equivariant. 

Conversely, a $P$--homomorphism $J^r(E)_o\to F_o$ extends uniquely to
a $G$--homo\-morphism $J^r(E)\to F$ and thus gives rise to an
invariant differential operator. Thus, invariant differential
operators $\Ga(E)\to\Ga(F)$ of order $\leq r$ are in bijective
correspondence with $P$--homomorphisms $J^r(E)_o\to F_o$. To avoid the
restriction on the order, one can simply pass to infinite jets and we
obtain that invariant differential operators $\Ga(E)\to\Ga(F)$ are in
bijective correspondence with $P$--homomorphisms $J^\infty(E)_o\to
F_o$, which factorize over some $J^r(E)$.

Surprisingly, the problem of determining all such homomorphisms has a
nice reformulation in term of (infinite--dimensional) representation
theory, which has led to a complete solution in several cases. Namely,
suppose that $E$ and $F$ correspond to representations $\Bbb E$ and
$\Bbb F$ of $P$, respectively. For the dual representation $\Bbb E^*$,
one can form the induced module $\Cal U(\frak g)\otimes_{\Cal U(\frak
p)}\Bbb E^*$, which is a $(\frak g,P)$--module, i.e.\ it admits
compatible actions of $\frak g$ and $P$. In the case where $\frak
p\subset\frak g$ is the Borel subalgebra (i.e.\ the minimal parabolic)
and $\Bbb E$ is irreducible, these are the {\em Verma--modules\/}
while for general $\frak p$ and irreducible $\Bbb E$, they are called
{\em generalized Verma--modules\/}. By a dualization argument and
Frobenius reciprocity one shows that for $\Bbb E$ and $\Bbb F$
irreducible, the space of all $P$--module homomorphisms
$J^\infty(E)_o\to F_o$, which factorize over some $J^r(E)_o$ is
isomorphic to the space of all $(\frak g,P)$--homomorphisms $\Cal
U(\frak g)\otimes_{\Cal U(\frak p)}\Bbb F^*\to\Cal U(\frak
g)\otimes_{\Cal U(\frak p)}\Bbb E^*$. Since these considerations are
essential for understanding of the links of our development to the
standard Bernstein--Gelfand--Gelfand resolutions, we 
provide some more details in Appendix \ref{A}.

Let us remark however that while there is a complete classification of
homomorphisms of Verma--modules in the complex case in \cite{BGG}, the
classification of homomorphisms of generalized Verma modules is a very
difficult problem, which is unsolved in general (even in the complex
case). There is a complete classification in the case of real rank one
for one dimensional representations in \cite{Lepowsky} and for general
representations in \cite{BoeColl1} and \cite{BoeColl2}. The problem in
the case of generalized Verma modules is the following: One has a
class of homomorphisms which are induced by homomorphisms of Verma
modules, the so--called standard homomorphisms. These are exactly the
homomorphisms which occur in Bernstein--Gelfand--Gelfand
resolutions. But it may happen that a homomorphism of Verma modules
induces the zero--homomorphism between generalized Verma modules, and
in this situation there may still be nonzero homomorphisms (the so
called non--standard homomorphisms). 

\subsection{Parabolic geometries}\label{2.8}
Some geometries can be viewed as curved analogs of the
homogeneous spaces $G/P$ considered above. For the purpose of this
paper, the best way to define them is simply as generalized spaces in
the sense of E.~Cartan. 

Let $\frak g=\frak g_{-k}\oplus\dots\oplus\frak g_k$ be a
real $|k|$--graded Lie algebra and let $G$ be a Lie group
with Lie algebra $\frak g$. Let $G_0$ and $P$ be the 
subgroups of $G$ defined in \ref{2.3} above. Then we define a (real)
{\em parabolic geometry of type $(G,P)$} on a smooth manifold $M$ to be a
principal $P$--bundle $\Cal G\to M$ equipped with a Cartan connection
of type $(G,P)$,
i.e.\ a differential form $\om\in\Om^1(\Cal G,\frak g)$ such that 
\begin{enumerate}
\item[(1)] $\om(\ze_X)=X$ for all $X\in \frak p$
\item[(2)] $(r^b)^*\om=\text{Ad}(b^{-1})\o \om$ for all $b\in P$
\item[(3)] $\om|_{T_u\Cal G}: T_u\Cal G\to \frak g$ is a linear
isomorphism for all $u\in \Cal G$.
\end{enumerate}

Here $\ze_X$ denotes the fundamental vector field generated by $X\in
\frak p$ and $r^b$ denotes the principal right action of $b\in
P$. Thus, $\om$ gives a smooth $P$--equivariant trivialization of the tangent
bundle of $\Cal G$, which reproduces the generators of fundamental
fields. Each $X\in {\frak g}$ defines the \emph{constant
vector field} $\om^{-1}(X)$ given by $\om^{-1}(X)(u)=\om_u^{-1}(X)\in T_u\Cal
G$. Clearly, a parabolic geometry of type $(G,P)$ on $M$ can only
exist if $M$ has the same dimension as $G/P$. 

In the complex setting, the Lie algebras and groups, as well as the manifold
$M$ are complex and the above definition remains unchanged except for
the replacement of {\em smooth} by {\em holomorphic}. Thus a {\em complex
parabolic geometry of the type $(G,P)$} on a complex manifold $M$ is given
by a holomorphic principal fiber bundle equipped with a holomorphic absolute
parallelism $\om$ with the three properties listed above.

The (real or complex) homogeneous space $G/P$ always carries a canonical
parabolic geometry, namely $\Cal G=G$ and the Cartan connection is given by
the left Maurer Cartan form. Then the constant vector fields
are exactly the left invariant fields on $G$.

It is fairly easy to make the parabolic geometries as defined above
into a category. Let $(\Cal G,\om)$ be a real parabolic geometry on $M$ and
$(\Cal G',\om')$ be a parabolic geometry on $M'$, and suppose that
$\Ph:\Cal G\to\Cal G'$ is a smooth homomorphism of principal $P$--bundles,
such that the induced map $\underline\Ph:M\to M'$ is a local
diffeomorphism. Then for any point $u\in\Cal G$ the tangent map
$T_u\Ph:T_u\Cal G\to T_{\Ph(u)}\Cal G'$ is a linear isomorphism, and
using this, one immediately verifies that $\Ph^*\om':=\om'\o T\Ph$ is
a Cartan connection on $\Cal G$. Now we define a {\em morphism\/} from
$(\Cal G,\om)$ to $(\Cal G',\om')$ to be a homomorphism $\Ph$ of
principal bundles such that the induced map $\underline\Ph:M\to M'$ is
a local diffeomorphism and such that $\om=\Ph^*\om'$. For complex parabolic 
geometries we additionally require all maps to be holomorphic.

Note that any homomorphism $\Ph:\Cal G\to\Cal G'$ of principal bundles
which lies over a local diffeomorphism can be viewed as a morphism 
$(\Cal G,\Ph^*\om')\to (\Cal G',\om')$. More generally, if $(\Cal
G',\om')$ is a parabolic geometry on $M'$ and $f:M\to M'$ is a local
diffeomorphism, then we can form the pullback bundle $f^*\Cal G'\to
M$. Then there is an induced homomorphism $\Ph:f^*\Cal G'\to\Cal G'$
of principal bundles which lies over $f$, and we get an induced
morphism $(f^*\Cal G',\Ph^*\om')\to (\Cal G',\om')$. 

\subsection{}\label{2.9}
For some purposes, the category of parabolic geometries as defined
above is too large, and one has to impose certain
restrictions. Usually, these restrictions are on the curvature of the
Cartan connection. Initially, the curvature of a Cartan connection
$\om$ is defined as the $\frak g$--valued two--form $K\in\Om^2(\Cal
G,\frak g)$ defined by the structure equation
$$K(\xi,\eta)=d\om
(\xi,\eta)+[\om(\xi),\om(\eta)],
$$ 
where $\xi$ and $\eta$ are vector
fields on $\Cal G$ and the bracket is in $\frak g$. Using the
properties of $\om$ one immediately verifies that $K$ is horizontal
and equivariant. In particular, this implies that $K$ is uniquely
determined by the {\em curvature--function\/} $\ka:\Cal G\to
\La^2\frak g_-^*\otimes \frak g$ defined by
$\ka(u)(X,Y)=K(u)(\om_u^{-1}(X),\om_u^{-1}(Y))$. There are two natural
ways to split $\ka$ into components. First, the splitting of $\frak g$
induces a splitting of $\ka$ according to the values in $\frak g$. In
particular, we can split $\ka=\ka_-\oplus\ka_{\frak p}$ according to
the splitting $\frak g=\frak g_-\oplus\frak p$. Following the classical 
terminology for affine connections, $\ka_-$ is called the {\em torsion\/} 
of $\om$. The other possibility is to
split $\ka$ according to homogeneous components. We denote the
homogeneous component of degree $i$ of $\ka$ by $\ka^{(i)}$. So
$\ka^{(i)}$ maps $\frak g_j\otimes\frak g_k$ to $\frak g_{i+j+k}$. 

Another important point is that the space $\La^2\frak
g_-^*\otimes\frak g$ is the second chain group $C^2(\frak g_-,\frak
g)$ in the standard complex for the Lie algebra cohomology $H^*(\frak
g_-,\frak g)$ of the nilpotent Lie algebra $\frak g_-$ with
coefficients in the ${\frak g_-}$--module $\frak g$. 
As we shall recall in detail in
Section \ref{4}, there is the adjoint $\partial^*$ to
the Lie algebra differential $\partial$ in this complex, so in
particular, we have $\partial^*:\La^2\frak g_-^*\otimes\frak g\to\frak
g_-^*\otimes\frak g$. 

\subsection*{Definition} 
Let $(\Cal G,\om)$ be a (real or complex) parabolic geometry 
on a manifold $M$, and let
$\ka$ be the curvature of $\om$. Then the parabolic geometry is called 
\begin{enumerate}
\item[(1)] {\em normal\/} if 
$\partial^*\o\ka=0$.
\item[(2)] {\em regular\/} if it is normal
and $\ka^{(i)}=0$ for all $i\leq 0$.
\item[(3)] {\em torsion--free\/} if 
$\ka_-=0$.
\item[(4)] {\em flat\/} if $\ka=0$.
\end{enumerate}

Note that forming the curvature of a Cartan connection is a natural
operation. This means that if $\Ph:\Cal G\to\Cal G'$ is a homomorphism
of principal bundles and $\om'$ is a Cartan connection with curvature
$K'$ and curvature--function $\ka'$ then the curvature $K$ and
curvature function $\ka$ of the pullback $\Ph^*\om'$ are given by
$K=\Ph^*K'$ and $\ka=\ka'\o\Ph$, respectively. Since all the
subclasses of parabolic geometries defined above are given by
restricting the values of the curvature--function, morphisms into a
parabolic geometry from one of the four subclasses can only come from
geometries from the same subclass. Clearly, for any of the four
subclasses the geometries belonging to the class form a full
subcategory of the category of all parabolic geometries of fixed type.  

\subsection{Examples}\label{2.10}
Before we review the construction of parabolic geometries
from underlying data, we present two well known examples.

\subsection*{Conformal structures}
Consider $\Bbb R^n$ with coordinates $x_1,\dots,x_n$ and the standard
inner product $\langle\ ,\ \rangle$ of signature $(p,q)$, and $\Bbb
R^{n+2}$ with coordinates $x_0,x_1,\dots,x_n,x_\infty$ and the inner
product associated to the quadratic form
$2x_0x_\infty+\langle(x_1,\dots,x_n),(x_1,\dots,x_n)\rangle$, which has
signature $(p+1,q+1)$. Let $G=SO_0(p+1,q+1)$ be 
the connected component of the special orthogonal group of this
metric. Then the Lie algebra $\frak g$ of $G$ admits a $|1|$--grading
by decomposing matrices into blocks of sizes $1$, $n$, and $1$, see
e.g. \cite[3.3(2)]{CSS1}. The construction of the canonical Cartan
connection for manifolds endowed with a conformal structure of
signature $(p,q)$, originally due to E.\ Cartan (see
\cite{Cartan-conf}), shows that conformal structures of this signature
are precisely the same thing as normal parabolic geometries
corresponding to that choice of $G$ and $P$. See \cite{CSS1} for a
construction of the canonical Cartan connection on conformal manifolds
in a style similar to the approach of this paper. In this special
situation, normal Cartan connections turn out to be automatically
regular and torsion free, so three of the four subclasses defined in
\ref{2.9} above coincide. The flat parabolic geometries in this case
are exactly the locally conformally flat manifolds. 

\subsection*{Partially integrable almost CR--structures}
The complex analog of the above construction leads to the
partially integrable almost CR--structures which present another
example of real parabolic geometries. Here we have to consider the
complex vector space  $\Bbb C^n$ with the
standard Hermitian inner product of signature $(p,q)$ and $\Bbb
C^{n+2}$ with the Hermitian inner product associated to $z_0\bar
z_\infty+\bar
z_0z_\infty+\langle(z_1,\dots,z_n),(z_1,\dots,z_n)\rangle$. Now we put
$G=PSU(p+1,q+1)$ the quotient of the special unitary group
corresponding to this Hermitian inner product by its center. Splitting
the matrices in the Lie algebra $\frak g$ of $G$ into blocks of sizes
$1$, $n$, and $1$ this time gives rise to a $|2|$--grading. The
construction of canonical Cartan connections in \cite{CS} shows that
partially integrable almost CR--structures with non--degenerate
Levi--form of signature $(p,q)$ are exactly the same thing as regular
parabolic geometries corresponding to $G$ (see \cite[4.14]{CS}). In
this case, three of the four subclasses of geometries defined in
\ref{2.9} above are really different: The torsion free parabolic
geometries in this case are precisely the CR--structures (see
\cite[4.16]{CS}), and the flat ones are those which are locally
isomorphic to the homogeneous model. The only coincidence in this case
is that normal parabolic geometries are automatically regular. 

\subsection{Underlying structures}\label{2.13} 
These two examples already show that identifying a
geometrical structure on a manifold as a parabolic geometry should be
rather the result of a theorem than a definition. In fact one can
show in a fairly general setting that certain parabolic geometries are
determined by underlying structures. This is the subject of the paper
\cite{CS} which generalizes \cite{Tan}, see also \cite{Mor} and
\cite{Yam}. To review the results, we first describe the underlying
structures we have in mind. 

Suppose that $(\Cal G,\om)$ is a regular parabolic geometry on a manifold
$M$. The first thing we get out of this is a filtration
$TM=T^{-k}M\supset T^{-k+1}M\supset\dots\supset T^{-1}M$ of the
tangent bundle of $M$. This is given by defining $T^iM$ to be the set
of those tangent vectors $\xi$ on $M$ for which there is a tangent vector
$\tilde\xi$ in $T\Cal G$ lying over $\xi$ with $\om(\tilde\xi)\in\frak
g_i\oplus\dots\oplus\frak g_k$. The latter condition is independent of 
the choice of
$\tilde\xi$ since changing the vector with fixed footpoint adds a
vertical vector whose image under $\om$ lies in $\frak p$, while
changing the footpoint leads to the adjoint action of an element of
$P$, which by definition preserves the subspace $\frak
g_i\oplus\dots\oplus\frak g_k$. Clearly, this filtration has the
property that the rank of $T^iM/T^{i+1}M$ equals the dimension of
$\frak g_i$ for all $i=-k,\dots,-1$. 

Now the underlying structures basically are given by considering the
bundles $\Cal G/P_+^i\to M$ for $i=1,\dots,k$ and the ``traces'' of
the Cartan connection that remain on these bundles. This ``trace'' on
the bundle $\Cal G/P_+^i\to M$ is a frame form of length $i$ in the
sense of \cite[3.2]{CS}. For the case $i=1$ the geometric meaning of
such a frame form is particularly easy to describe: It is exactly a
reduction to the structure group $G_0$ of the associated graded vector
bundle 
$$
\operatorname{gr}TM=T^{-k}M/T^{-k+1}M\oplus\dots\oplus T^{-2}M/T^{-1}M\oplus T^{-1}M
$$ 
to the tangent bundle $TM$. The fact that the curvature--function
$\ka$ of the regular Cartan connection $\om$ has the property that
$\ka^{(i)}=0$ for all $i\leq 0$ is reflected in a property of the
underlying frame forms called the structure equation, see
\cite[3.4]{CS}. The bundle $\Cal G/P_+^i$ together with the frame form
of length $i$, which satisfies the structure equations is called the
underlying $P$--frame bundle of degree $i$. Again, for $i=1$ this
condition can be easily understood geometrically. It is equivalent
to the fact that the algebraic Lie bracket on $\operatorname{gr}TM$ which
comes from the reduction to the group $G_0$ is induced by the Lie
bracket of vector fields, that is it is given by a (generalized) Levi--form. 

Now the main result of \cite{CS} can be stated (with the help of the
language of Dynkin diagrams for the pairs $(\frak g,\frak p)$
mentioned in \ref{2.2} above) as follows: 

{\em Let $(\frak g,\frak p)$, $G$, $P$, and $G_0$  be as in \ref{2.3} and
suppose throughout that no simple factor of $\frak g$ is contained in $\frak
g_0$ and $\frak g$ does not contain a simple factor of type
$A_1$. Then:\newline 
(1) If $(\frak g,\frak p)$ does not contain any simple factor of one of the 
types
$$
\begin{picture}(65,10)
\put(0,2.5){\line(1,0){20}}
\put(60,2.5){\line(-1,0){20}}
\put(0,2.5){\makebox(0,0){$\times$}}
\put(15,2.5){\makebox(0,0){$\bullet$}}
\put(30,2.5){\makebox(0,0){$\cdots$}}
\put(45,2.5){\makebox(0,0){$\bullet$}}
\put(60,2.5){\makebox(0,0){$\bullet$}}
\end{picture} 
\quad\text{or}\quad 
\begin{picture}(65,10)
\put(0,2.5){\line(1,0){20}}
\put(60,3.5){\line(-1,0){15}}
\put(60,1.5){\line(-1,0){15}}
\put(45,2.5){\line(-1,0){5}}
\put(0,2.5){\makebox(0,0){$\times$}}
\put(15,2.5){\makebox(0,0){$\bullet$}}
\put(30,2.5){\makebox(0,0){$\cdots$}}
\put(45,2.5){\makebox(0,0){$\bullet$}}
\put(60,2.5){\makebox(0,0){$\bullet$}}
\put(52.5,2.5){\makebox(0,0){\smaller[3]$\langle$}}
\end{picture}
$$
then any regular parabolic geometry can be
reconstructed from the underlying $P$--frame bundle of degree one, and
any $P$--frame bundle of degree one comes from a regular parabolic
geometry. Thus, in all these cases regular parabolic geometries are
the same thing as manifolds with filtered tangent bundle plus
reductions of $\operatorname{gr}TM$ to the group $G_0$ such that the resulting
algebraic bracket is induced by the Lie bracket.\newline
(2) If $\frak g$ contains simple factors of one of the two above
types, then any regular parabolic geometry can be reconstructed from
the underlying $P$--frame bundle of degree two and any such bundle
comes from a regular parabolic geometry. Moreover, any $P$--frame
bundle of degree one can be extended (in various ways) to a $P$--frame
bundle of degree two.} 

The classical examples of the second case are the projective
structures where the $P$--frame bundle of degree one is simply the
full frame bundle and all the structure is contained in the choice of
an extension to a $P$--frame bundle of degree two.
The other exceptional examples are the so called {\em projective contact
structures}.

\subsection{Natural bundles and operators}\label{2.11}
We will not go into much detail in the generalities about natural
bundles and natural operators, but just outline the basic facts. We do
not want to compare the various notions of naturality (this will be
taken up elsewhere) but just show that the operators we are going to
construct are natural (or invariant) in any reasonable sense.

Given a representation of $P$ on a vector space $\Bbb V$ and a
parabolic geometry $(\Cal G\to M,\om)$ we can form the associated
bundle $VM=\Cal G\x_P\Bbb V\to M$. If $\Ph:\Cal G\to\Cal G'$ is a
homomorphism of principal bundles which covers a local diffeomorphism
$\underline\Ph:M\to M'$, then we get an induced homomorphism of vector
bundles $VM\to VM'$ which lies over the same map $\underline\Ph$ and
restricts to a linear isomorphism in each fiber. To put it in another
way, we get a functor from the category of parabolic geometries to the
category of vector bundles over manifolds of the same dimension as
$G/P$ and vector bundle homomorphisms which cover local
diffeomorphisms and induce linear isomorphisms in each fiber such that
the composition of the base functor with the given functor equals the
base functor. Thus, we get a special case of a gauge natural bundle as
defined in \cite[Chapter XII]{KMS}. 

Consider next a fixed category of real parabolic geometries, and two
representations $\Bbb V$ and $\Bbb W$ of $P$. Let $V$ and $W$ be the
corresponding natural vector bundles. A {\em natural linear operator\/}
mapping sections of $V$ to sections of $W$ is defined to be a system
of linear operators $D_{(\Cal G,\om)}:\Ga(VM)\to\Ga(WM)$, where
$M$ is the base of $\Cal G$ such that for any morphism $\Ph:(\Cal
G,\om)\to (\Cal G',\om')$ we have 
$$
\Ph^*\o D_{(\Cal G',\om')}=D_{(\Cal G,\om)}\o\Ph^*.
$$ 
This definition implies immediately, that each of the
operators is local both in the section and in the Cartan connection:
Suppose that $s\in\Ga(VM)$ vanishes identically on an open subset
$U\subset M$. Then there is an obvious inclusion morphism $i:(\Cal
G|_U,\om|_U)\to (\Cal G,\om)$ and $i^*s=0$. Thus also $i^*(D_{(\Cal
G,\om)}(s))=0$, i.e.\ {}$D_{(\Cal G,\om)}(s)$ is identically zero on
$U$. Similarly, assume that $\om$ and $\om'$ are two Cartan
connections which coincide on $\Cal G|_U$. Then for any section
$s\in\Ga(VM)$ we have $D_{(\Cal G,\om)}(s)|_U=D_{(\Cal
G,\om')}(s)|_U$. In particular, the classical Peetre theorem implies
that each of the operators $D_{(\Cal G,\om)}$ is locally over $M$ a
finite order differential operator with respect to the arguments in
the vector bundles and the Cartan connection.  

For complex parabolic geometries, we deal with holomorphic representations of
$P$, the natural vector bundles are holomorphic, and the natural operators
act on holomorphic sections. Let us also remark that all these
concepts extend to non-linear objects without essential changes.

\subsection{}\label{2.12}
The natural operators on the category of flat parabolic geometries are
particularly easy to describe: It is a classical result on Cartan
connections that any flat parabolic geometry is locally isomorphic to
the homogeneous model $G/P$ (see \cite[4.12]{CS} for a proof in the
setting of parabolic geometries). This immediately implies that any
natural operator on the category of flat parabolic geometries is
uniquely determined by its value on the homogeneous model
$G/P$, i.e. the parabolic geometry $(G\to G/P, \om)$. 
Moreover, an operator on the flat model extends to a natural
operator on the category of flat parabolic geometries if and only if
it is natural with respect to all automorphisms of $(G,\om)$. The left
multiplication by any element of $G$ induces an automorphism of the
principal bundle $G\to G/P$ and by left invariance of the Maurer
Cartan form this actually is an automorphism of the parabolic geometry
$(G,\om)$. On the other hand, by \cite[Theorem 3.5.2]{Sharpe} the only
smooth functions $G\to G$ which pull back the Maurer Cartan form to
itself are the constant left translations. Thus $G$ is exactly the
group of all automorphisms of $(G,\om)$. But this immediately implies
that an operator on the homogeneous model extends to a natural
operator on the category of flat parabolic geometries if and only if
it is invariant in the sense of definition \ref{2.5}. Thus for the
flat case, the description of natural operators is equivalent to a
problem in representation theory.

Usually, the question on more general natural operators is then posed
(in the special cases that have been studied so far) as the question
of the existence of {\em curved analogs} of invariant operators. This should
be viewed as follows: As we discussed in \ref{2.6}, an invariant
operator of order $r$ is induced by a $P$--module 
homomorphism $J^r(E)_o\to
F_o$, which does not factor over $J^{r-1}(E)_o$. Now the kernel of the
projection $J^r(E)_o\to J^{r-1}(E)_o$ is the bundle
$S^rT^*(G/P)\otimes E$, so it corresponds to the representation
$S^r\frak p_+\otimes \Bbb E$. Thus the invariant operator gives rise
to a $P$--module homomorphism $S^r\frak p_+\otimes\Bbb E\to \Bbb F$,
which in turn gives a $G$--equivariant homomorphism between the
corresponding homogeneous  vector bundles
which is precisely the symbol of the operator we started with. But
this $P$--module homomorphism induces a homomorphism of associated bundles on
any parabolic geometry, so for any parabolic geometry $(\Cal G,\om)$
over a manifold $M$, we get the corresponding homomorphism
$S^rT^*M\otimes EM\to FM$. Now a curved analog of an invariant
operator is a natural operator such that for each $(\Cal G,\om)$ the
symbol of $D_{(\Cal G,\om)}$ is the above homomorphism. Otherwise put,
the question is whether we can extend a given natural operator from
the category of flat parabolic geometries to some larger category of
parabolic geometries without changing its symbol, which, as a natural
transformation, makes sense on any parabolic geometry. 

\subsection{}\label{2.14} We conclude this introductory section with some
more remarks on the beautiful geometric structure underlying each 
parabolic geometry. This topic deserves much more attention than we could
pay here and it will be studied in detail elsewhere. Some first steps
have been done in \cite{Sl97}.  

Suppose that $(\Cal G,\om)$ is a real parabolic geometry on a manifold
$M$. Then we have the tower of principal fiber bundles 
$\Cal G\to\Cal G/P_+\to M$ and the top level has the structure group 
$P_+$. Now using the Baker--Campbell--Hausdorff formula, Proposition
\ref{2.3} can be restated in the form that for any $g\in P$ there
is a unique $g_0\in G_0$ and a unique $Z\in\frak p_+$ such that
$g=g_0\exp(Z)$. But using this, one easily shows that the bundle $\Cal
G\to\Cal G/P_+$ admits global $G_0$--equivariant smooth sections.  Namely,
one can use a local trivialization of $\Cal G\to M$ to construct equivariant
sections over the preimage in $\Cal G/P_+$ of appropriate open subsets of
$M$. Such local sections can then be glued to a global section using a
partition of unity (compare with the proof of \cite[Lemma 3.6]{CSS1}). As in
this last reference one also proves that the space of all these sections is
an affine space modeled on the space
$\Om^1(M)$ of one--forms on $M$. 

$$
\xymatrix@R=3mm@C=20mm{
{\Cal G} \ar[r] 
&{\Cal G/\Cal P_+} \ar[r] \ar@(ul,ur)[l]^{\si}
& M
\\
{\om\strut} \ar@{.>}[u]
&{\si^*(\om_{{\frak g}_-}+\om_{{\frak g}_0})}
\ar@{.>}[u]
}$$
Each such global section $\si$ reduces the structure group of the tangent space
$TM$ to $G_0$ and induces an affine connection $\ga^\si=
\si^*(\om_{{\frak g}_-}+\om_{{\frak g}_0})$ on $TM$. This affine connection
is $\si$--related to another Cartan connection $\om^{\si}$ on $\Cal G$, which
differs from $\om$ only in the ${\frak p}_+$--component. 
The class of all connections $\ga^\si$ is a
straightforward generalization of Weyl structures on conformal
geometries and all differential operators built of the Cartan connection
$\om$ can be expressed by uniform formulae in terms of these affine
connections and their torsions and curvatures. The technique based on
this general framework was developed 
systematically for all $|1|$-graded algebras $\frak g$ in
\cite{CSS1,CSS2,CSS3}. 

\section{Semi--holonomic jet modules and strongly invariant
operators}\label{3} 

Semi--holonomic jet prolongations of modules were first introduced in
the context of AHS--structures in \cite{CSS1}. Here we develop the
concept in the more general setting of parabolic geometries and we discuss
how the homomorphisms of semi--holonomic jet prolongations give rise to
natural operators. Throughout this section, there will be essentially no
differences in the arguments for the real and complex parabolic geometries.
Thus we shall not mention the field of scalars explicitly, and one has to
think of the proper real or complex modules in the applications below. 

\subsection{The absolutely invariant derivative}\label{3.1}
Suppose that $(\Cal G,\om)$ is a parabolic geometry on a manifold
$M$. We mentioned in \ref{2.5}, that the tangent and cotangent bundles on the
homogeneous spaces are homogeneous vector bundles. The Cartan connection
$\om$ extends this identification to all parabolic geometries as follows:

We identify $\frak g_-$ (as a $P$--module) with $\frak g/\frak p$, and
consider the map $\Cal G\x \frak g_-\to TM$ defined by mapping $(u,X)$
to $Tp\cdot\om_u^{-1}(X)$, where $p:\Cal G\to M$ is the
projection. The equivariancy of the Cartan connection immediately
implies that this factors to a vector bundle homomorphism $\Cal
G\x_P\frak g_-\to TM$. Since this is immediately seen to be
surjective, it must be an isomorphism of vector bundles by dimensional
reasons. Thus we have identified $TM$ with the natural bundle associated to
the  $P$--module ${\frak g_-}$.
Now, the invariance of the Killing form on $\frak g$ implies that $\frak
g/\frak p$ and $\frak p_+$ with the actions induced by the adjoint
action are dual $P$--modules. Thus, similarly as above the cotangent
bundle $T^*M$ of $M$ can be identified with the bundle $\Cal
G\x_P\frak p_+$ (implicitly, this has been used in \ref{2.14} above). 

There is a nice way to encode the action of vector fields on functions
(or equivalently the exterior derivative of functions) using the
identifications made above. As we have seen, a typical tangent
vector on $M$ can be written as $Tp\cdot\om_u^{-1}(X)$ for an element
$X\in\frak g_-$. Acting with this tangent vector on a smooth function
$f\in C^\infty(M,\Bbb R)$, we get $\om_u^{-1}(X)\cdot (f\o p)$. Now, 
smooth functions on $M$ are in
bijective correspondence with smooth $P$--invariant functions on
$\Cal G$, the correspondence given by mapping $f$ to $f\o p$. 
To any smooth, $P$--invariant function $f$ on $\Cal G$ we
associate a function $\nabla^\om f:\Cal G\to L(\frak g_-,\Bbb
R)$ defined by $\nabla^\om f(u)(X):=\om_u^{-1}(X)\cdot f$. The
equivariancy properties of $\om$ imply that the map $\nabla^\om f$ is
$P$--equivariant. Taking into account the above identification of $T^*M$
with an associated bundle and of $L(\frak g_-,\Bbb R)\simeq\frak
p_+$, we see that $\nabla^\om f$ is a one form on $M$, which
by definition coincides with $df$.

The above procedure immediately suggests a generalization. Let $\Bbb
V$ be any representation of $P$ and let $VM=\Cal G\x_P\Bbb V$ be the
corresponding associated bundle. Then we can identify smooth sections
of $VM$ with smooth maps $\Cal G\to\Bbb V$, which are
$P$--equivariant. Now to any smooth function $s:\Cal G\to\Bbb V$ we
associate a smooth function $\nabla^\om s:\Cal G\to L(\frak
g_-,\Bbb V)$ defined by 
$$\nabla^\om s(u)(X):=\om_u^{-1}(X)\cdot s.
$$ 
Obviously, this defines a differential operator 
$$C^\infty(\Cal G,\Bbb V)\to C^\infty(\Cal G,L(\frak g_-,\Bbb V))
$$ and these
operators (for all $(\Cal G,\om)$) form a natural operator on all
parabolic geometries in the sense of \ref{2.11}. This operation is
called the {\em universal covariant derivative\/} in the book
\cite[p. 194]{Sharpe}. In \cite[2.3]{CSS1} we have chosen to call it the
{\em absolutely invariant derivative\/}. The reason for the latter
name also shows the main drawback of this operation: It is not really
covariant, i.e. if one starts with an equivariant map $s$ (i.e. a
section of $VM$) the result is not equivariant in general. Thus in general, 
if we start with a section, the result of the invariant derivative is not a
section of a bundle anymore. 

\subsection{}\label{3.2}
There is a way, however, to make a section of an associated bundle out
of a section of an associated bundle and its absolutely invariant
derivative. This is called the {\em invariant one--jet\/} of the
section. To describe it, we first have to analyze the action of
$G$ on one--jets in the homogeneous case. Thus, let us consider a
representation $\Bbb V$ of $P$, the corresponding homogeneous bundle
$V(G/P)=G\x_P\Bbb V$ and its first jet prolongation $J^1(V(G/P))\to G/P$.  
As we noted in \ref{2.6} this is again a homogeneous bundle, and
we want to describe the corresponding action of $P$ on its standard
fiber $\Cal J^1(\Bbb V):=J^1(V(G/P))_o$. As we noticed in \ref{2.4}
it suffices to understand this space as a module over $G_0$ and over
$\frak p_+$ (in fact, already $\frak g_1$ would be sufficient).

If we think of sections in $\Ga(V(G/P))$ as $P$--equivariant
functions $s\in C^\infty(G,\Bbb V)^P$,  then the 1--jets of sections at
the distinguished point $o\in G/P$ are identified with 1--jets of
equivariant functions at the unit $e\in G$ and the action is
given by $g.(j^1_e s)= j^1_e(s\o\ell_{g^{-1}})$ for all $g\in G$. 
Thus, the induced action of $Z\in\frak p$ on the section $s$ is given
by the differentiation in the direction of the right invariant vector
field $R_Z$ on $G$, $Z.j^1_e s=-j^1_e(R_Z\cdot s)$.

Now we can identify a one--jet $j^1_e(s)$ with $(s(e),ds(e))$ and as
we saw in \ref{3.1} above, $ds(e)=\nabla^\om s(e)$. As a vector
space we can thus write 
$$\Cal J^1(\Bbb V)=\Bbb V\oplus(\frak g_-^*\otimes \Bbb V)
$$ 
and we have to understand the induced actions
of $G_0$ and $\frak p_+$ on this space. Let us first assume that $g\in
G_0$. Then $(s\o\ell_{g^{-1}})(e)=s(g^{-1})=g\cdot s(e)$ by equivariancy
of $s$. On the other hand, we have to evaluate $\om_e^{-1}(X)\cdot
(s\o\ell_g^{-1})$. This can be computed as 
\begin{multline*}
\ddt s(g^{-1}\exp(tX))=\ddt s(g^{-1}\exp(tX)gg^{-1})=\\
=\om_e^{-1}(\Ad(g^{-1})X)\cdot (g\cdot
s)=g\cdot(\om_e^{-1}(\Ad(g^{-1})X)\cdot s) .
\end{multline*}
Now since $g\in G_0$, we have $\Ad(g^{-1})X\in\frak g_-$ for all
$X\in\frak g_-$ (the adjoint action on $\frak g_-$ coincides with the
induced action on $\frak g/\frak p$ in this case), so we see that 
$\Cal J^1(\Bbb V)=\Bbb V\oplus(\frak g_-^*\otimes \Bbb V)$ even as a
$G_0$--module. 

For $Z\in\frak p_+$ we have $-(R_Z\cdot s)(e)=Z\cdot (s(e))$ by
the infinitesimal version of equivariancy of $s$. On the other hand, for the
derivative component we have to compute the linear mapping 
${\frak g_-}\ni X\mapsto-\om^{-1}(X)\cdot R_Z\cdot s(e)$. Since
$\om^{-1}(X)$ is left invariant, it commutes with $R_Z$ and the
resulting expression depends only on $R_Z(e)=Z=\om^{-1}(Z)(e)$, and we
get 
\begin{align*}
-\om^{-1}(X)\cdot R_Z\cdot s(e)&=-\om^{-1}(Z)\cdot\om^{-1}(X)\cdot
s(e)\\
&=-\om^{-1}(X)\cdot\om^{-1}(Z)\cdot
s(e)-[\om^{-1}(Z),\om^{-1}(X)]\cdot s(e).  
\end{align*}
The infinitesimal version of equivariancy of $s$ shows that the first
term in the last expression gives $Z\cdot (\om_e^{-1}(X)\cdot
s(e))$. Since $\om^{-1}(\_)$ is just the left invariant vector field,
the second term gives $-\om_e^{-1}([Z,X])\cdot s$. Now let us split
$\ad(Z)=\ad_-(Z)\oplus \ad_\frak p(Z)$ according to the splitting
$\frak g=\frak g_-\oplus\frak p$. Then the $\ad_\frak p(Z)(X)$--part
acts algebraically by equivariancy of $s$ while the rest simply
produces $-\om_e^{-1}(\ad_-(Z)(X))\cdot s$. 

Thus, if we denote elements of $\Cal J^1(\Bbb V)$ as pairs $(v,\ph)$,
where $v\in\Bbb V$ and $\ph$ is a linear map from $\frak g_-$ to $\Bbb
V$, then the appropriate action of $Z\in\frak p_+$ is given by 
$$
Z\cdot (v,\ph)=(Z\cdot v, X\mapsto Z\cdot (\ph(X))-\ph(\ad_-(Z)(X))+
\ad_\frak p(Z)(X)\cdot v),
$$
i.e.\ we get the tensorial action plus one additional term mapping the
value--part to the derivative--part. 

This action can also be nicely written
in a tensorial notation. To do this let us choose a basis
$\{\eta_\al\}$ of $\frak p_+$ such that each element $\eta_\al$ is
homogeneous of degree $|\eta_\al|$, and let $\{\xi_\al\}$ be the dual
basis of 
$\frak g_-$ (with respect to the Killing form $B$). Now consider an
element $(v_0,Z_1\otimes v_1)\in \Cal J^1(\Bbb V)$, where
$v_0,v_1\in\Bbb V$ and $Z_1\in\frak p_+\cong\frak g_-^*$. 
Then by definition $Z_1\otimes v_1$ maps $X\in\frak g_-$ to
$B(Z,X)v_1$. Thus $[Z,X]_-:=\ad_-(Z)(X)$ is mapped to
$B(Z_1,[Z,X]_-)v_1$. Since the Killing form vanishes on $\frak
p_+\x\frak p$, this can be rewritten as
$B(Z_1,[Z,X])v_1=B([Z_1,Z],X)v_1$. 
Moreover, we can write $\ad_Z$ as an element of $L(\frak
g_-,\frak g)\cong\frak p_+\otimes\frak g$ in the form
$\sum_\al\eta_\al\otimes [Z,\xi_\al]$. This implies that for $Z$
homogeneous of degree $|Z|$, we may rewrite the action on $\Cal J^1{\Bbb V}$
as
$$
Z\cdot (v_0,Z_1\otimes v_1)=(Z\cdot v_0,Z_1\otimes Z\cdot
v_1+[Z,Z_1]\otimes v_1+\sum_{|\eta_\al|\leq |Z|}\eta_\al\otimes
[Z,\xi_\al]\cdot v_0). 
$$ 

A simple computation shows that $\Cal J^1(\_)$ can be made into a
functor on the category of $P$--modules by defining 
$$
\Cal J^1(f)(v,\ph):=(f(v),f\o\ph)
$$ 
for each $P$--module homomorphism $f:\Bbb V\to \Bbb W$.

\subsection{}\label{3.3}
Surprisingly, the first jet prolongation of representations introduced above
leads for any parabolic geometry to a natural identification of the first jet
prolongation of any natural bundle with an associated bundle, i.e. with
another natural bundle.  Let $(\Cal
G,\om)$ be a parabolic geometry on $M$, let $\Bbb V$ be a representation of
$P$, and let $VM$ be the corresponding associated bundle over $M$.

\begin{prop*} 
The invariant differential $\nabla^\om$ defines the mapping
$$
\io:C^{\infty}(\Cal G,\Bbb V)^{P}\to C^\infty(\Cal G, \Cal J^1\Bbb V)^P, \quad
\io(s)(u)= (s(u), (X\mapsto \nabla^\om s(u)(X)))
$$
which yields an isomorphism $J^1VM\simeq \Cal G\x_P \Cal J^1\Bbb V$.

For each fiber bundle map $VM\to WM$ induced by a $P$--module
homomorphism $f:\Bbb V\to \Bbb W$, the first jet prolongation of the
bundle map is induced by the $P$--module homomorphism $\Cal J^1(f)$. 
\end{prop*}
\begin{proof} Let us recall that $\nabla^\om s(u)(X)=
\om^{-1}(X)(u)\cdot s$. Thus the mapping $\io: s\mapsto (s,\nabla^\om
s)$ is well defined and depends on first jets only, so we only have to
check that the values are actually equivariant. First, for $g\in G_0$
we have to compute $(s(u\cdot g),\nabla^\om s(u\cdot
g))$. Equivariancy of $s$ implies $s(u\cdot g)=g^{-1}\cdot(s(u))$. The
second component maps $X\in\frak g_-$ to $\om_{u\cdot g}^{-1}(X)\cdot
s$. Now the equivariancy of $\om$ immediately implies that
$\om_{u\cdot g}^{-1}(X)=Tr^g\cdot \om_u^{-1}(\Ad(g)X)$. Since $g\in
G_0$ we see that $\Ad(g)X\in\frak g_-$ and using equivariancy of $s$
again, we see that $\nabla^\om s(u\cdot g)$ maps $X$ to $g^{-1}\cdot
(\om_u^{-1}(\Ad(g)X)\cdot s)$, and thus $(s(u\cdot g),\nabla^\om s(u\cdot
g))=g^{-1}\cdot (s(u),\nabla^\om s(u))$. 

On the other hand, we have to check equivariancy for the infinitesimal
action of $Z\in\frak p_+$. Thus, we have to compute $((\ze_Z\cdot
s)(u),\ze_Z\cdot (\nabla^\om s)(u))$. Equivariancy of $s$ implies that
the first component equals $-Z\cdot (s(u))$. The second component maps
$X\in\frak g_-$ to $(\ze_Z\cdot \om^{-1}(X)\cdot s)(u)$. Now
$\ze_Z=\om^{-1}(Z)$ and we can rewrite the expression as
$$
(\om^{-1}(X)\cdot \om^{-1}(Z)\cdot
s)(u)+[\om^{-1}(Z),\om^{-1}(X)]\cdot s(u).
$$
Since the curvature of $\om$ is horizontal and $\om^{-1}(Z)$ is
vertical, we may rewrite the second term in this expression as
$(\om^{-1}([Z,X])\cdot s)(u)$. Now we can split $[Z,X]$ into a $\frak
g_-$ and a $\frak p$--component and conclude as in \ref{3.2} above
that $((\ze_Z\cdot
s)(u),\ze_Z\cdot (\nabla^\om s)(u))=-Z\cdot (s(u),\nabla^\om s(u))$.

Clearly, this construction gives a smooth injective homomorphism of
vector bundles $J^1VM\to \Cal G\x_{P}\Cal J^1\Bbb V$, which covers the
identity map on $M$. Since both bundles clearly have the same rank,
this must be an isomorphism. 

Finally, consider a homomorphism $f:\Bbb V\to \Bbb W$. The
corresponding bundle map $VM\to WM$ is induced by $(u,v)\mapsto
(u,f(v))$, and so the induced action on sections is induced by  
$$
s \mapsto (x\mapsto (u(x),f\o s(u(x)))).
$$
Taking 1--jet of this expression we obtain just the homomorphism $\Cal
J^1(f)$. 
\end{proof}

\subsection{Semi--holonomic jets}\label{3.4}
Since we posed no conditions on the representation $\Bbb V$ above, we
can iterate the functors $J^1$ on the associated vector bundles as
well as the functors $\Cal J^1$ on the $P$--modules. Proposition
\ref{3.3} then implies that the $r$--th iteration $J^1\dots J^1 VM$ is
an associated bundle to $\Cal G$ corresponding to the  $P$--module 
$\Cal J^1\dots \Cal J^1\Bbb V$. Let us look more carefully
at $\Cal J^1\Cal J^1\Bbb V$ and $J^1J^1VM$. There are two obvious
$P$--module homomorphisms $\Cal J^1\Cal J^1 \Bbb V\to \Cal J^1\Bbb V$,
the first one given by the projection $p_{\Cal J^1\Bbb V}$ defined on
each first jet prolongation by projection to the first component, and
the other one obtained by the action of $\Cal J^1$ on $p_{\Bbb V}$.
Thus there is the submodule $\bar\Cal J^2\Bbb V$ in $\Cal J^1\Cal
J^1\Bbb V$ on which these two projections coincide. As a vector space
and a $G_0$--module we have 
$$
\bar\Cal J^2\Bbb V=\Bbb V\oplus (\frak g_-^*\otimes \Bbb V)\oplus
(\frak g_-^*\otimes \frak g_-^*\otimes \Bbb V).
$$
The two $P$--module homomorphisms $\Cal J^1(p_{\Bbb V})$ and $p_{\Cal
J^1\Bbb V}$ give rise to vector bundle homomorphisms $J^1J^1VM\to
J^1VM$ which are just the two standard projections on the second
non--holonomic jet prolongation. So we conclude that the second
semi--holonomic prolongation $\bar J^2VM$ is naturally isomorphic to 
$\Cal G\x_P \bar\Cal J^2\Bbb V$.

Iterating this procedure, we obtain the $r$--th semi--holonomic jet
prolongations and $\Cal J^1(\bar\Cal J^{r}\Bbb V)$ equipped with two
natural projections onto $\Cal J^1(\bar\Cal J^{r-1}\Bbb V)$, 
which correspond to
the usual projections on the first jet prolongation of semi--holonomic
jets. Their equalizer is then the submodule $\bar\Cal J^{r+1}\Bbb
V$. As a $G_0$--module 
$$
\bar\Cal J^r\Bbb V= \bigoplus_{i=0}^r(\otimes^i\frak g_-^*\otimes \Bbb
V).
$$

\begin{prop*} 
For each positive integer $r$, the $r$--th
semi--holonomic jet prolongation $\bar J^rVM$ carries the natural
structure of associated vector bundle $\Cal G\x_P \bar\Cal J^r\Bbb
V$.  Moreover, there is the natural embedding
\begin{gather*}
J^rVM\to \bar J^rVM\simeq \Cal G\x_P \bar\Cal J^r\Bbb V\\
j^rs(u)\mapsto \{u,(s(u), \nabla^{\om}s(u),\dots, (\nabla^{\om})^rs(u))\}
.
\end{gather*}
\end{prop*}
\begin{proof} The first part of the statement has been already shown.
What remains is to discuss the equivariancy properties of the invariant
differentials. However also this follows from the first order case easily by
induction, using only the definition of the semi--holonomic prolongations.
\end{proof}

\subsection{Strongly invariant operators}\label{3.5}
The problem, why we cannot work with true (holonomic) $r$--jets but
have to use the semi--holonomic ones, is that absolutely invariant
derivatives commute only for flat Cartan connections. More
precisely, from the definition of the absolutely invariant derivative
and the properties of the curvature, one immediately concludes the so
called {\em general Ricci--identity\/}
\begin{align*}
(\nabla^\om\nabla^\om s)(u)(X\otimes Y-Y\otimes X)&=
\nabla^\om s(u)([X,Y])+\ka_{\frak p}(X,Y)\cdot(s(u))\\
&\quad-\nabla^\om s(u)(\ka_-(X,Y))
\end{align*}
for all $X,Y\in {\frak g_-}$.
This also shows that the torsion--part of $\ka$ has a quite different
geometric meaning than the component valued in $\frak p$.  
Thus, the identification from proposition \ref{3.4} has values in the
$P$--submodule $\Cal J^r(\Bbb V)$ of symmetric elements
$\oplus_{i=0}^r(S^i\frak g_-^*\otimes \Bbb V)$ in the flat case. 
Consequently we have recovered the standard 
identification of the $r$--th holonomic
jet prolongation of a homogeneous bundle with an associated
bundle for flat geometries, but this does not work in the curved case. 

Nevertheless, one can well use the semi--holonomic jet prolongations
to generate invariant operators. Suppose that $\Bbb V$ and $\Bbb W$
are representations of $P$ and suppose that $\Ph:\bar\Cal J^r(\Bbb
V)\to\Bbb W$ is a homomorphism of $P$--modules. Then for any parabolic
geometry $(\Cal G,\om)$ we can define a differential operator
$\Ga(VM)\to\Ga(WM)$ as follows: For a section $s$ viewed as an
equivariant function $\Cal G\to\Bbb V$ define 
$$
D_{(\Cal G,\om)}(s)(u)=\Ph(s(u),\nabla^\om
s(u),\dots,(\nabla^\om)^rs(u)).
$$
{From} Proposition \ref{3.4} above it follows that this gives a section
of the bundle $WM$ and that each $D_{(\Cal G,\om)}$ is a differential
operator of order $\leq r$. Moreover, by construction the operators
$D_{(\Cal G,\om)}$ form a natural operator on the category of all
parabolic geometries in the sense of \ref{2.11}. Operators arising in
this way will be called {\em strongly invariant operators\/} in the
sequel. We will often not distinguish carefully between a strongly
invariant operator and the corresponding homomorphism $\bar\Cal
J^r(\Bbb V)\to\Bbb W$. Thus, the semi--holonomic jet modules give a
possibility to construct natural operators for a parabolic geometry in
a completely algebraic way, since one only has to construct a
homomorphism between two finite dimensional $P$--modules. 

There is a slight problem about strongly invariant operators,
however. Namely, even if a homomorphism $\bar\Cal J^r(\Bbb V)\to\Bbb
W$ does not factor over $\bar\Cal J^{r-1}(\Bbb V)$, the corresponding
operators may be of order strictly less than $r$ or even identically
zero. To see this, note that we can easily compute the symbol of a
strongly invariant operator. This symbol is a vector bundle
homomorphism $S^rT^*M\otimes VM\to WM$, which is induced by a
homomorphism $S^r\frak g_-^*\otimes\Bbb V\to\Bbb W$. Using Proposition
\ref{3.4} it is clear that this homomorphism is given by restricting
$\Ph$ to $S^r\frak g_-^*\otimes\Bbb V$, viewed as a submodule of
$\otimes^r\frak g_-^*\otimes\Bbb V$, which in turn can be viewed as a
submodule of $\bar\Cal J^r(\Bbb V)$. Thus, if a homomorphism restricts
to zero on the symmetric part of the top component of the jet--module,
then the corresponding operator actually is of lower order (and
contains terms involving the curvature of the Cartan connection). 

There is an important situation in which this problem does not play
any role. Suppose that we have an operator of order $r$ in the flat
case with nontrivial symbol, and suppose that we can find a
homomorphism $\bar\Cal J^r(\Bbb V)\to\Bbb W$ which induces this
operator (in the flat case). Then this gives a curved analog of the
operator in question, and there is no problem with the symbol at
all. This will always be the case for the operators we are going to
study. In particular, since $\bar\Cal J^1(\Bbb V)=\Cal J^1(\Bbb V)$,
any first order invariant operator on the category of flat parabolic
geometries is automatically strongly invariant, and thus has a
canonical curved analog.  

\begin{rem}\label{3.6}\em
There are operators which are natural (invariant) in the sense of
\ref{2.11} but are not strongly
invariant. Basically, there is only one example of such an operator
known: It is shown in \cite{GJMS} that on conformal manifolds of
dimension $2m$ there exists a conformally invariant 
$m$--th power of the Laplacian on smooth functions. In
\cite{ES} it is shown that this operator is not strongly
invariant. It can, however, be written in terms of absolutely
invariant derivatives, and thus it is also natural. 
In fact, it is shown in \cite{Sl96} that for
AHS--structures, i.e.\ parabolic geometries corresponding to
$|1|$--graded Lie algebras, naturality of (even non-linear) operators 
is equivalent to the possibility to express them by means of the absolute 
invariant derivative and curvature of the defining Cartan connection, 
and this, in turn, is equivalent to the existence of a
universal formula in terms of all underlying affine connections, cf.
\ref{2.14}.

The existence of invariant operators which are not strongly invariant
is due to symmetries of the curvature of a Cartan connection. Suppose
that we write an expression in terms of absolutely invariant
derivatives and check whether the result is
$P$--equivariant. Otherwise put, we can compute the obstruction
against being equivariant which usually contains expressions involving
the curvature of the Cartan connection and its derivatives. In the
case of a strongly invariant operator, these obstructions vanish
algebraically. But the jets of the curvature of any Cartan connection
have certain symmetries, basically due to the Bianchi identity, see
e.g. \cite[4.9]{CS}. This implies that expressions that do not vanish
algebraically, still may vanish whenever the jet of the curvature of a
Cartan connection is inserted, and this is precisely what happens in
the case of the critical powers of the Laplacian.
\end{rem}

\subsection{Twisted invariant operators}\label{3.8}
Besides the completely reducible representations (which come from the
reductive subgroup $G_0$) there is a second class of particularly
simple representations of the group $P$. Namely one can take a
representation of the full (semisimple) group $G$ and restrict it to
$P$. These representations have particularly nice features in the case of
the flat model since they give rise to trivial homogeneous
bundles. There are many ways to see that, but the most appropriate one
for our purposes is to associate to any element $v$ in a
representation $\Bbb V$ of $G$ a global nonzero section of the
associated bundle $G\x_P\Bbb V$. To do this, we just have to specify
a $P$--equivariant map $G\to\Bbb V$, and we define this map simply by
$g\mapsto g^{-1}\cdot v$. This map is even $G$--equivariant and not
only $P$--equivariant.  

There is a simple generalization of this result. Suppose
that $\Bbb W$ is any representation of $P$. Then sections of
$W(G/P)$ are in bijective correspondence with $P$--equivariant maps
$G\to\Bbb W$. Now we define a map on sections of homogeneous bundles
\begin{gather*}
\Ga (W(G/P))\otimes \Bbb V\to\Ga\bigl(W(G/P)\otimes V(G/P)\bigr)\\
s\otimes v \mapsto (g\in G \mapsto s(g)\otimes g^{-1}\cdot v)
\end{gather*}
and one immediately
verifies that this is an isomorphism of $G$--modules. In particular,
this implies that if $\Bbb W'$ is another $P$--representation and
$D:\Ga(W(G/P))\to \Ga(W'(G/P))$ is an invariant differential operator,
then we can pull back 
$$
D\otimes\id_\Bbb V:\Ga(W(G/P))\otimes\Bbb V\to
\Ga(W'(G/P))\otimes\Bbb V
$$
along these isomorphisms to get an
invariant operator 
$$
D_\Bbb V:\Ga\bigl(W(G/P)\otimes V(G/P)\bigr)\to\Ga\bigl(W'(G/P)\otimes
V(G/P)\bigr).
$$
This operator is called the twisted invariant
operator corresponding to $D$ and $\Bbb V$. 

Now, let us notice that the above isomorphism between the spaces of sections
of the associated bundles induces a $P$--module isomorphism $\bar\Cal
J^r(\Bbb W)\otimes {\Bbb V}\simeq \bar\Cal J^r({\Bbb W}\otimes{\Bbb V})$ for
all $P$--modules ${\Bbb W}$ and $G$--modules ${\Bbb V}$ and all orders $r$.
Thus, for strongly invariant operators $D$, 
we may extend the construction of the twisted invariant operators to natural
operators $D_{\Bbb V}$ acting on all geometries $(\Cal G,\om)$ of the type
$(G,P)$ and the resulting operators are
again strongly invariant. Let us remark that a completely algebraic
treatment of this construction has been worked out (in the special case 
of the AHS-structures) in \cite{Cap}.

In particular, we obtain the strongly invariant twisted operators $D_{\Bbb
V}$ for all first order invariant operators $D$ on the homogeneous vector
bundles and all $G$--modules ${\Bbb V}$.

\subsection{Twisted exterior derivatives}\label{3.8a}
The standard exterior derivatives $d$ on the differential forms on 
$G/P$ are first order invariant
operators (since they are even invariant under the action of all
diffeomorphisms of $G/P$), so we can apply the construction above to get the
twisted exterior derivatives 
$$
d_\Bbb V:\Ga\bigl(\La^nT^*(G/P)\otimes
V(G/P)\bigr)\to\Ga\bigl(\La^{n+1}T^*(G/P)\otimes V(G/P)\bigr)
$$
for $n=0,\dots,\text{dim}(G/P)$. Moreover, the operators $d_{\Bbb V}$
are strongly invariant, since they are of first order, and so there
are the corresponding $P$--module homomorphisms on 
the semi--holonomic jet modules. Since we will need it later, 
we will compute these homomorphisms explicitly. 

Let us start
with the ordinary exterior derivative. We have already noted in
\ref{3.1} that the exterior derivative of functions equals the
absolutely invariant derivative. To compute the exterior derivative
for general differential forms, we first have to describe nicely the
identification of $n$--forms with smooth equivariant functions
$G\to\La^n\frak p_+$. Throughout, we are going to identify $\La^n\frak
p_+$ with the space of $n$--linear alternating maps from $\frak
g_-\cong \frak g/\frak p$ to $\Bbb K$. Now using the identification of
the tangent bundle of $G/P$ with $G\x_P\frak g_-$ described in
\ref{3.1}, one easily verifies that the relation between a form
$\ph\in\Om^k(G/P)$ and the corresponding function $s:G\to\La^n\frak
p_+$ is given by 
$$
(p^*\ph)(g)(\om_g^{-1}(X_1),\dots,\om_g^{-1}(X_n))=s(g)(X_1,\dots,X_n),
$$
where $p^*\ph$ is the pullback of $\ph$ along the projection $p:G\to
G/P$, and the $X_i$ are in $\frak g_-$. Note that this formula remains
correct for $X_i\in\frak g$ if one interprets $s(g)$ as an $n$--linear
map on $\frak g$ which vanishes if at least one argument lies in
$\frak p$. 

\begin{lem*}
Let $s$ and $ds$ be the functions on $G$ corresponding to differential forms 
$\ph$ and $d\ph$ on $G/P$, respectively. Then the formula for the exterior
derivative reads as
\begin{align*}
ds(X_0,\dots,X_n)=\ 
&\sum_{i=0}^n(-1)^i(\nabla^\om s)(g)(X_i)(X_0,\dots,\hat i,\dots,X_n)+\\
&\sum_{i<j}(-1)^{i+j}s(g)([X_i,X_j],X_0,\dots,\hat i,\dots,\hat
j,\dots,X_n)
\end{align*}
where $\om$ is the left Maurer-Cartan form on $G$ and, as usual, 
the hat denotes omission. 
\end{lem*}
\begin{proof}
To compute the function corresponding to $d\ph$, we just have to
evaluate $p^*(d\ph)(g)=d(p^*\ph)(g)$ on vector fields of the form
$\tilde X(g)=\om_g^{-1}(X)$. We have
\begin{align*}
d(p^*\ph)&(\tilde X_0,\dots,\tilde X_n)=\sum_{i=0}^n(-1)^i\tilde
X_i\cdot (p^*\ph)(\tilde X_0,\dots,\hat i,\dots\tilde X_n)+\\
&+\sum_{i<j}(-1)^{i+j}(p^*\ph)([\tilde X_i,\tilde X_j],\tilde
X_0,\dots,\hat i,\dots,\hat j,\dots,\tilde X_n).
\end{align*}
Inserting $p^*\ph$ from above
and evaluating at $g$, we see directly that the first summand agrees with
the first summand in the claimed formula,
which clearly equals $n+1$ times the alternation of $(\nabla^\om
s)(g)$ evaluated at $(X_0,\dots,X_n)$. 

For the second summand, we just
have to note that by the Maurer--Cartan equation for $\om$ we have
$[\tilde X_i,\tilde X_j]=\widetilde{[X_i,X_j]}$. Thus, this summand
gives exactly the other part of the required formula.
\end{proof}

Now let us pass to the general case of a $V(G/P)$--valued $n$--form,
where $\Bbb V$ is a representation of the whole group $G$. Any such
form can be written as a finite sum of expressions of the form
$\ph\otimes\tilde v$, where $\ph\in\Om^n(G/P)$ and $\tilde v$ is the
global section of $V(G/P)$ corresponding to $v\in \Bbb V$ as in
\ref{3.8} above. By definition, the twisted exterior derivative is
given by $d_\Bbb V(\ph\otimes\tilde v)=(d\ph)\otimes\tilde v$. Now let $s$
be the function corresponding to $\ph$ and denote by $\tilde v$ also
the function corresponding to the global section. From above, we thus
see that $d_\Bbb V(\ph\otimes\tilde v)$ is represented by the function
which maps $(X_0,\dots,X_n)$ to 
\begin{multline*}\tag{$*$}
\sum_{i=0}^n(-1)^i(\nabla^\om s)(g)(X_i)(X_0,\dots,\hat
i,\dots,X_n)\tilde v(g)+\\
+\sum_{i<j}(-1)^{i+j}s(g)([X_i,X_j],X_0,\dots,\hat i,\dots,\hat
j,\dots,X_n)\tilde v(g).
\end{multline*}
By definition of the absolutely invariant derivative, we have
$$
\nabla^\om(s\otimes\tilde v)(X)=\nabla^\om s(X)\otimes\tilde
v+s\otimes(\nabla^\om\tilde v(X))
$$
and the infinitesimal version of
$G$--invariance of $\tilde v$ says that 
$$\nabla^\om\tilde
v(g)(X)=-X\cdot (\tilde v(g)).
$$ 
Thus we may rewrite the first summand in \thetag{$*$} as 
\begin{multline}\tag{$**$}
\sum_{i=0}^n(-1)^i\nabla^\om(s\otimes\tilde v)(g)(X_i)(X_0,\dots,\hat
i,\dots,X_n)+\\
+\sum_{i=0}^n(-1)^iX_i\cdot (s(g)(X_0,\dots,\hat i,\dots,X_n)\tilde v(g)).
\end{multline}
Finally note that the second term in \thetag{$**$} adds up with the
second term in \thetag{$*$} to the value of the standard Lie algebra
differential $\partial:C^n(\frak g_-,\Bbb V)=\La^n\frak
g_-^*\otimes\Bbb V\to C^{n+1}(\frak g_-,\Bbb V)$ (cf. \ref{4.1} for the
explicit formula) applied to the map
$s(g)\otimes\tilde v(g)$ evaluated on $(X_0,\dots ,X_n)$. Thus we may
summarize:

\begin{prop}\label{3.8b}
The twisted exterior derivative $d_\Bbb V$ on $G/P$ is a strongly invariant
operator induced by the
$P$--module homomorphism
$\bar\Cal J^1(\La^n\frak p_+\otimes \Bbb V)\to \La^{n+1}\frak
p_+\otimes \Bbb V$, 
which is given by the formula
$$
(f_0,Z\otimes f_1)\mapsto\partial(f_0)+(n+1)Z\wedge f_1,
$$
where we view elements of $\La^n\frak p_+\otimes \Bbb V$ as
$n$--linear alternating maps from $\frak g_-$ to $\Bbb V$ and $Z\wedge
f_1$ denotes the alternation of the map 
$(X_0,\dots,X_n)\mapsto B(Z,X_0)f_1(X_1,\dots,X_n)$. 
\end{prop}

\begin{kor}\label{3.9}
The Lie algebra differential $\partial$ satisfies
$$
(W\cdot\partial(f)-\partial(W\cdot f))
=(n+1)\sum_{|\eta_\al|\le|W|}\eta_\al\wedge[W,\xi_\al]\cdot f
$$
for $f\in\La^n\frak p_+\otimes\Bbb V$ and $W\in\frak p_+$, where
$\xi_\al$ and $\eta_\al$ are homogeneous dual bases of $\frak g_-$ and
$\frak p_+$ with respect to the Killing form. 
\end{kor}
\begin{proof} The claim can be verified by a nice and elementary, but
tedious algebraic computation. However, the previous proposition
offers the following simple argument:

We know that the formula for the strongly invariant operator
$$
d_{\Bbb V}(f_0,Z\otimes f_1)=\partial(f_0)+(n+1)Z\wedge f_1
$$ 
is $P$--equivariant. Thus for all $f_0$, $f_1\in {\Bbb V}$, $Z\in
{\frak p}_+$, $W\in {\frak p}_+$ we obtain the equality of the
following two expressions 
\begin{align*}
d_{\Bbb V}(W\cdot(f_0,&Z\otimes f_1)) = d_{\Bbb V}
((W\cdot f_0,W\cdot(Z\otimes f_1) + \sum \eta_\al\otimes[W,\xi_\al]\cdot
f_0)=\\
&=\partial(W\cdot f_0)+(n+1)W\cdot(Z\wedge f_1) + 
(n+1)\sum \eta_\al\wedge[W,\xi_\al]\cdot f_0\\
W\cdot(\partial(f_0)&+(n+1)Z\wedge f_1)=W\cdot(\partial f_0) + 
(n+1)W\cdot(Z\wedge f_1).
\end{align*}
This yields the required formula.
\end{proof}

\subsection{The covariant exterior derivatives}\label{3.10}
Proposition \ref{3.8b} offers a canonical curved analog of the twisted
exterior derivatives on all manifolds with a parabolic
geometry of the type $(G,P)$.  It should be remarked that we may obtain
another curved analog as follows. For any parabolic geometry $(\Cal G,\om)$
on $M$, we consider the extended bundle $\tilde\Cal G=\Cal G\x_PG$, which is
a principal $G$--bundle over $M$. It is a classical observation that the
Cartan connection $\om$ induces a principal connection $\tilde\om$ on
$\tilde\Cal G$. Now if
$\Bbb V$ is a representation of $G$, then we can view the
corresponding natural bundle $VM=\Cal G\x_P\Bbb V$ also as $VM=\tilde\Cal
G\x_G\Bbb V$, and thus we have the induced linear connection on this
bundle. The covariant exterior derivative with respect to
this connection gives a natural operator on $VM$--valued forms on $M$. If
$s:\tilde {\Cal G}\to \La^k{\frak p}_+\otimes {\Bbb V}$ is the
equivariant function corresponding to 
a $k$-form $\ph$ on $M$, then the value of the 
latter operator is a $(k+1)$-form on $M$, given by the formula
\begin{align*}
d^{\tilde \om}s(u)(&X_0,\dots,X_n)=\sum_{i=0}^k(-1)^i
\nabla^{\tilde\om}_{X_i}s(u)(X_0,\dots,\hat i,\dots,X_k)+\\ 
&+\sum_{i<j}(-1)^{i+j}s(u)([X_i,X_j],X_0,\dots,\hat i,\dots,\hat
j,\dots,X_k) 
\end{align*}
where $X_0,\dots,X_k\in {\frak g}_-$, $u\in \tilde{\Cal G}$,
$\nabla^{\tilde \om}_{X_i}s(u)$ means the derivative of $s$ in the
direction of the horizontal vector at $u$ determined by $X_i$, and there are
the standard omissions of arguments in the expressions on the right hand side.
Indeed, $d^{\tilde\om}$ is defined as the pullback of the standard $d$ on
$\tilde{\Cal G}$ by the horizontal projection of $\tilde \om$, applied to the
pullback of the $k$-form $\ph$ on $M$ by the projection $p:\tilde {\Cal G}\to
M$. Since the curvature of $\tilde \om$ produces vertical fields on 
$\tilde{\Cal G}$, the above formula equals to the standard evaluation 
of $d(p^*\ph)$ on the horizontal lifts of vector fields on $M$.

These operators coincide with the twisted exterior derivatives on 
the homogeneous space but they differ in general. The explicit general 
comparison is as follows:

\begin{lem*}
Let ${\Bbb V}$ be a $G$-module, $VM$ the corresponding natural vector bundle
over a manifold $M$ equipped with a parabolic geometry $({\Cal G}, \om)$.
The covariant exterior derivative $d^{\tilde \om}$ on $\La^kT^*M\otimes VM$,
$k>0$,
and the twisted exterior derivative $d_{\Bbb V}$ on the same space satisfy
$$
d^{\tilde \om}\phi = d_{\Bbb V}\ph + i_{\ka_{-}}\ph$$
where $\ka_-$ is the torsion--component of the curvature of $\om$ and 
$i_{\ka_{-}}\ph$ is the usual insertion operator evaluated on $\ka_-$ and
$\ph$, i.e. the alternation of $\ph(\ka_{-}(X_0,X_1),X_2,\dots,X_k)$ over
the arguments.
\end{lem*}
\begin{proof}
The key to the required formula is in the expressions \thetag{$*$} and
\thetag{$**$} in \ref{3.8a}. Namely, the latter expressions which were derived
on the homogeneous spaces describe also the
twisted exterior derivatives in general, but we have to be aware that instead
of the bracket $[X_i,X_j]$ in \thetag{$*$} we have to plug in
$$
\om(u)([\om^{-1}(X_i),\om^{-1}(X_j)])= [X_i,X_j]-\ka(u)(X_i,X_j)
.$$
At the same time, for all $u\in{\Cal G}\subset \tilde{\Cal G}$, the
covariant derivative $\nabla^{\tilde\om}$ of a section $s:\tilde {\Cal G}\to
{\Bbb V}$ relates to the absolute invariant derivative as
$$
\nabla^{\tilde\om}s(u)(X)= \nabla^\om s(u)(X) + X\cdot s(u)
$$
(since the horizontal fields given by $\tilde \om$ equal to 
$\om^{-1}(X)$ minus the fundamental field $\zeta_X$).

Combining the latter two facts, we see that exactly the expression 
\begin{multline*}
i_{\ka_-}\ph(u)(X_0,\dots,X_k)=\\
\sum_{i<j}(-1)^{i+j}\ph(u)
(\ka_-(X_i,X_j),X_0,\dots,\hat i,\dots,\hat j,\dots,X_k)
\end{multline*}
has to be added to $d_{\Bbb V}(u)\ph(X_0,\dots,X_k)$ in order to obtain the
covariant derivative.  This is exactly the evaluation of the 
insertion operator, cf. \cite[8.2]{KMS}.
\end{proof}

The latter lemma shows that our twisted exterior differentials $d_{\Bbb V}$
are certain torsion adjusted versions of the standard covariant exterior
derivatives. In particular, even in the case ${\Bbb V}={\Bbb R}$ the twisted
derivative $d_{\Bbb R}$ equals to the usual exterior derivative
$d$ if and only if the geometry is torsion--free.

\subsection{Remarks}\label{3.11}
(1) As we saw in \ref{3.8a}, the isomorphism
$$
\Ga(W(G/P))\otimes\Bbb V\cong\Ga(W(G/P)\otimes V(G/P))
$$
of $G$--modules induces an isomorphism of
$P$--modules $\bar\Cal J^r(\Bbb W)\otimes\Bbb V\cong\bar\Cal J^r(\Bbb
W\otimes \Bbb V)$ for any $P$--representation $\Bbb W$ and
$G$--representation $\Bbb V$. This can also be proved algebraically
along the lines of \cite{Cap}. This isomorphism can then be used
to define twisted versions of any strongly invariant operators in a
completely algebraic way. Using
this picture, the subsequent developments in this paper can be viewed
as a curved analog of the Jantzen--Zuckermann translation principle in
representation theory. The first version of such a curved translation
procedure appeared in the context of 4--dimensional conformal geometry in
\cite{EasR}, see also \cite{EasConf}. 

(2) The twisted exterior derivatives give a sequence 
$$
\Ga(VM)\to\Om^1(M;VM)\to\dots\to\Om^{\text{max}}(M;VM)\to 0,
$$
of invariant differential operators, where sections and forms are
smooth in the real case and holomorphic in the complex case. In the
case of the flat model, this sequence is just the pullback of the
tensor product of the (smooth or holomorphic) de~Rham sequence with
$\Bbb V$, so it is a resolution of the constant sheaf $\Bbb V$. In the
case of a general parabolic geometry, it fails to be a
complex. Actually, it is easy to verify that the composition $d_{\Bbb
V}\o d_{\Bbb V}$ is just  given by the action of the curvature of
$\om$. Thus, in the case of a flat parabolic geometry, we still get a
complex, which by Lemma \ref{3.10} coincides with the complex given by
the covariant exterior derivative with respect to the flat linear
connection induced by the Cartan connection. Note however, that on a
flat parabolic geometry bundles corresponding to representation of $G$
are no more trivial in general.  

(3) As a $G_0$--module, one can split any representation $\Bbb W$ of
$P$ as $\oplus\Bbb W_j$ according to eigenvalues of the grading
element $E\in\frak g_0$. Clearly, the action of $\frak p_+$ maps
$\frak g_i\otimes\Bbb W_j$ to $\Bbb W_{j+i}$. In particular, we
can apply this to $\La^n\frak p_+\otimes\Bbb V$ to split the space
$\Om^n(M;VM)$ into homogeneous components, and analyze how the twisted
exterior derivative behaves with respect to this splitting. From the
formula in Proposition \ref{3.8b} it is obvious that $d_\Bbb V$ never
lowers homogeneous degree and the component of the same homogeneous
degree as the input is just the Lie algebra differential $\partial$
composed with the given form. Thus, the homogeneous component of
degree zero of $d_\Bbb V$ is algebraic and equals $\partial$. This
observation is crucial for the subsequent development. Using the fact
that the Lie algebra cohomology of $\frak g_-$ with coefficients in
$\frak g$ admits a Hodge theory (which we will discuss in the next
section), we will show that we can replace the sequence of remark (2)
above by a different sequence in which only sections of completely
reducible bundles occur, and which is a complex computing the same
cohomology if the original sequence was a complex. 

\section{Curved analogs of Bernstein--Gelfand--Gelfand resolutions}\label{4}
In this section, we first discuss the Hodge--structure on the standard
complex for the cohomology $H^*(\frak g_-,\Bbb V)$ for a $\frak
g$--module $\Bbb V$. Then we come to the core of the paper, the
construction of a huge class of distinguished natural operators on all
parabolic geometries.

\subsection{}\label{4.1}
We have already mentioned the standard complex for the cohomology
$H^*(\frak g_-,\Bbb V)$ in \ref{3.8a}. The chain groups in this complex
are the groups $C^n(\frak g_-,\Bbb V)=\La^n\frak g_-^*\otimes\Bbb V$,
which are viewed as the spaces of $n$--linear alternating maps from
$\frak g_-$ to $\Bbb V$. The differential 
$$
\partial:C^n(\frak g_-,\Bbb V)\to C^{n+1}(\frak g_-,\Bbb V)
$$ 
is defined by 
\begin{align*}
\partial(f)(X_0,\dots,X_n)=&\sum_{i=0}^n(-1)^i X_i\cdot
f(X_0,\dots,\hat i,\dots, X_n)+\\
+&\sum_{i<j}(-1)^{i+j}f([X_i,X_j],X_0,\dots,\hat i,\dots,\hat j,\dots,X_n),
\end{align*}
where the hats denote omission. Clearly, if we start with a
representation $\Bbb V$ of the group $G$, then $\partial$ is a
homomorphism of $G_0$--modules, and it is well known that
$\partial\o\partial=0$. 

The crucial fact for us is that on this standard complex there is a
Hodge theory, which was first introduced for complex simple Lie
algebras in \cite{Kostant}. The most conceptual way to describe this Hodge
structure is to use the natural duality between $\frak g_-$ and $\frak
p_+$ via the Killing form. This is a duality of $G_0$--modules, but if
we consider $\frak g_-$ as a $P$--module via the adjoint action and
the identification with $\frak g/\frak p$, then it even is a duality
of $P$--modules by invariance of the Killing form. Thus, given a
representation $\Bbb V$ of $\frak g$ and its dual $\Bbb V^*$, we can
naturally identify $C^n(\frak p_+,\Bbb V^*)$ with the dual $P$--module
of $C^n(\frak g_-,\Bbb V)$. Thus, the dual map to the Lie algebra
differential $\partial:C^n(\frak p_+,\Bbb V^*)\to C^{n+1}(\frak
p_+,\Bbb V^*)$ can be viewed as a linear map 
$$
\partial^*:C^{n+1}(\frak g_-,\Bbb V)\to C^n(\frak g_-,\Bbb V)
$$ 
which is called the {\em
codifferential\/}. From the above, it is obvious that the
codifferential is a $G_0$--homomorphism and
$\partial^*\o\partial^*=0$. Moreover, one immediately verifies that
the Lie algebra differential for $\frak p_+$ is even a
$P$--homomorphism and thus the same is true for $\partial^*$.

A formula for $\partial^*$ can be easily computed for elements of the
form $Z_0\wedge\dots\wedge Z_n\otimes v$, where the $Z_i$ are in
$\frak p_+$ and $v$ is in $\Bbb V$. Pairing this element with a
multilinear map $\ps\in C^{n+1}(\frak p_+,\Bbb V^*)$, we simply get
$\ps(Z_0,\dots,Z_n)(v)$. Using this, one immediately computes that 
\begin{align*}
\partial^*(Z_0\wedge\dots\wedge Z_n\otimes v)&=\sum_{i=0}^n(-1)^{i+1} 
Z_0\wedge\cdots\hat i\dots\wedge Z_n\otimes Z_i\cdot v+\\
&+\sum_{i<j}(-1)^{i+j}[Z_i,Z_j]\wedge\cdots\hat i\cdots\hat j\dots\wedge
Z_n\otimes v.
\end{align*}
{}From this formula, it is again obvious that $\partial^*$ is a
$P$--homomorphism. 

Using Lie theory, one constructs an inner product on the spaces of
cochains, with respect to which $\partial$ and $\partial^*$ are
adjoint operators. The proof for this fact in the generality we
need it is only a rather simple extension of results available in the
literature, see e.g. \cite{Tan, Yam}. For the sake of completeness and
the convenience of the reader, we give a complete proof in Appendix
\ref{B}.  

\subsection{}\label{4.4}
This adjointness result has a number of important consequences: First
of all one gets a harmonic theory for the cohomology $H^*(\frak
g_-,\Bbb V)$. We define the Laplacian
$$
\square=\partial\o\partial^*+\partial^*\o\partial.
$$ 
Then for each $n$
this is a $G_0$--endomorphism of $C^n(\frak g_-,\Bbb V)$. Moreover, the
adjointness implies that
$\ker(\square)=\ker(\partial)\cap\ker(\partial^*)$ and we have a
$G_0$--invariant splitting 
$$C^n(\frak g_{-},V)=\im(\partial)\oplus
\ker(\square)\oplus\im(\partial^*).
$$ 
This implies then that the
cohomology group $H^n(\frak g_-,\Bbb V)$ is isomorphic (as a $G_0$--module)
to the subspace $\ker(\square)\subset C^n(\frak g_{-},\Bbb
V)$. Moreover, the situation between $\partial$ and $\partial^*$ is
completely symmetric, so we can as well compute the cohomology groups
$H^*(\frak g_{-},\Bbb V)$ as $\ker(\partial^*)/\im(\partial^*)$. This is
more suitable for our purposes, since, as we have noticed above,
$\partial^*$ is even a $P$--homomorphism. This also implies that (even
as a $G_0$--module) the cohomology group $H^n(\frak g_{-},{\Bbb V})$ is dual
to $H^n(\frak p_+,\Bbb V^*)$. 

Thus, we get a canonical action of $P$ on the cohomology groups
$H^n(\frak g_-,\Bbb V)$. We claim, that this module is completely
reducible, i.e.\ a direct sum of irreducibles. To prove this, we only
have to show that $\frak p_+$ acts trivially on the cohomology
groups. Fortunately, there is the following simple observation 

\begin{lem*} Let $Z\in\frak p_+$ and $f\in
C^n(\frak g_-,\Bbb V)\cong \La^n\frak p_+\otimes\Bbb V$. Consider
$Z\cdot f\in C^n(\frak g_-,\Bbb V)$ and $Z\wedge f\in C^{n+1}(\frak
g_-,\Bbb V)$, as defined in \ref{3.8b}. Then
$$
\partial^*(Z\wedge f)=-Z\cdot f - Z\wedge\partial^*(f).
$$
\end{lem*}
\begin{proof} This is a direct consequence of the general formula for
$\partial^*$ in \ref{4.1}.
\end{proof}
Now, one immediately concludes that the $\frak p_+$--action
maps $\ker(\partial^*)$ to $\im(\partial^*)$, and thus in particular
the induced action on the cohomology groups is trivial.

In \cite{Kostant}, B.~Kostant computed the cohomology groups
$H^*(\frak p_+,\Bbb V)$ in the case when $\frak g$ is complex and
simple and $\Bbb V$ is a complex irreducible representation. 
The basic idea in the proof is to analyze the action of the Laplacian
$\square$ in terms of Casimir operators. 

In fact, our construction of the sequences of natural operators will not
need the explicit knowledge of the cohomologies. 
On the other hand, the full
power of Kostant's theorem is necessary in order to compare the resulting 
sequences with
the standard BGG--resolutions in representation theory. 

Let us also remark here, that the knowledge of the second cohomologies
with values in $\frak g$ determines nicely the structure of the
curvature of normal parabolic geometries, see e.g. \cite{Yam,SlSchm}. 

\subsection{A sketch of the construction}\label{5.1}
Let us return to the twisted de~Rham sequence 
$$
\Ga(VM)\to\Om^1(M;VM)\to\dots\to\Om^{\text{max}}(M;VM)\to 0
$$
on a manifold $M$ equipped with a parabolic geometry $(\Cal G,\om)$.
For each $i$, $\Om^i(M;VM)$ is the space of sections (smooth in the real
case, holomorphic in the complex setting) of the natural
bundle associated to the representation $\La^i\frak p_+\otimes\Bbb
V$. The maps $\partial$, $\partial^*$, and $\square$ defined above 
induce maps on the spaces of sections, which we denote by
the same symbols. Moreover, since these maps are induced by pointwise
operations the Hodge decomposition of $\La^i\frak p_+\otimes\Bbb V$
gives rise to a Hodge decomposition
$$
\Om^i(M;VM)=\im(\partial)\oplus\ker(\square)\oplus\im(\partial^*).
$$
One has to be careful, however, that this decomposition is not
$P$--invariant but just $G_0$--invariant, since $\partial^*$ is a
$P$--homomorphism but $\partial$ and $\square$ are not. Thus the latter
decomposition makes explicit geometrical sense only after a reduction of
$\Cal G$ to $G_0$, cf. the discussion in \ref{2.14}.

Since $\partial^*$ is a $P$--homomorphism, the kernel
$\ker(\partial^*)$ and the image $\im(\partial^*)$ are the spaces of
sections of natural subbundles of $\La^nT^*M\otimes VM$. Moreover,
from \ref{4.4} we know that the quotient
$\ker(\partial^*)/\im(\partial^*)$ can be identified with the space of
sections of the bundle associated to the (completely reducible)
representation $\Bbb H^n_\Bbb V=H^n(\frak g_-,\Bbb V)$ of $P$, so we
get an algebraic natural operator from the subset $\ker(\partial^*)$
of $\Om^n(M;VM)$ to the space of smooth sections of the natural bundle
corresponding to the representation $\Bbb H^n_\Bbb V$. If $\Bbb E$ is
an irreducible component of $\Bbb H^n_\Bbb V$, then we can further
project onto this component to get an algebraic natural operator
$\ker(\partial^*)\to\Ga(EM)$. 

On the other hand, $\Bbb H^n_\Bbb V$ can be
identified with $\ker(\square)\subset \La^n\frak p_+\otimes\Bbb V$ 
as a $G_0$--module, so
we may view any section of the corresponding bundle as a $VM$--valued
$n$--form, but this is not a natural operator. The main point of the
following will be that one can construct a natural differential
operator $L$ from sections of the bundle corresponding to $\Bbb H^n_\Bbb
V$ to $VM$--valued $n$--forms in $\ker(\partial^*)$, 
which has this inclusion as the lowest homogeneous component. 
Otherwise put, one can split the algebraic
projection $\pi$ constructed above by a natural differential operator $L$. 
Moreover, it will turn out that this operator is fully determined by the
following surprising fact: {\it For each section $\al\in \Gamma(H^n_{\Bbb
V}M)$ there exists the unique section $L(\al)\in
\operatorname{ker}(\partial^*)\subset \Om^n(M;VM)$ such that  $\pi\o
L(\al)=\al$ and $d_{\Bbb V}(L(\al))\in
\operatorname{ker}(\partial^*)\subset\Om^{n+1}(M;VM)$.}
\begin{equation*}
\xymatrix{
& {\ker(\partial^*)} \ar@<-1.5ex>[d]_{\pi}
& {\ker(\partial^*)} \ar@<-1.5ex>[d]_{\pi} & {\ }
\\
{\dots}\ar@{.>}[r] \ar[ru]^-{d_\Bbb V\o L}
& {\Gamma(H^i_{\Bbb V}M)} \ar@<-1.5ex>[u]_L \ar[ru]^-{d_\Bbb
V\o L}\ar@{.>}[r]
& {\Gamma(H^{i+1}_{\Bbb V}M)} \ar@<-1.5ex>[u]_L \ar@{.>}[r]
&{\dots}}
\end{equation*}
Summarizing the prospective achievement,
the twisted exterior derivatives will produce plenty of natural 
differential operators indicated by the dotted arrows in the diagram.

The idea for the construction of this natural differential operator $L$ is
fairly simple. Recall from \ref{3.8b} that the lowest homogeneous
component of $d_\Bbb V$ equals the Lie algebra differential
$\partial$. Suppose we have a section $s$ in the bundle corresponding to
$\Bbb H^n_\Bbb V$, which is homogeneous of some degree $i$. Then it lies
in $\ker(\square)$ and thus in particular in $\ker(\partial)$, so the
homogeneous component of degree $i$ of  $d_\Bbb V(s)$ is automatically
zero. The idea is now to extend $s$ to $\tilde s$ in
such a way that $d_\Bbb V(\tilde s)$ is as small as possible. The
homogeneous component of degree $i+1$ of $d_\Bbb V(s)$ can be split
into components in $\im(\partial)$, $\ker(\square)$, and
$\im(\partial^*)$, and the best we can do to kill it is to add to $s$
an element $s_{i+1}$ which is homogeneous of degree $i+1$ such that
$\partial(s_{i+1})$ is the negative of the $\im(\partial)$--component
of $d_\Bbb V(s)$ in degree $i+1$. There is a freedom in the choice of
$s_{i+1}$ which can be fixed by requiring that
$s_{i+1}\in\im(\partial^*)$ (which is a complement to $\ker(\partial)$
by the adjointness). But this allows us already to compute $s_{i+1}$:
Since $\partial^*(s_{i+1})=0$ we see that $\square
(s_{i+1})=\partial^*(\partial(s_{i+1}))$. But $\partial(s_{i+1})$ is
just the negative of the $\im(\partial)$--part of the homogeneous
component of degree $i+1$ of $d_\Bbb V(s)$, so this is
known. Moreover, by definition $\square$ commutes both with $\partial$
and $\partial^*$, and $\ker(\square)\cap \im(\partial^*)=\{0\}$. Thus
$\square$ restricts to an isomorphism
$\im(\partial^*)\to\im(\partial^*)$. Hence we can compute $s_{i+1}$ by
applying $\square^{-1}\o\partial^*$ to the homogeneous component of
degree $i+1$ of $d_\Bbb V(s)$. Similarly we can continue to add an
appropriate homogeneous component of degree $i+2$ and so on. 

{}From this description it is {\em not at all obvious} that this construction
produces a natural operator, since the map $\square^{-1}$ involved in
the construction is not a $P$--homomorphism, and the subsequent steps
of the construction use $d_\Bbb V-\partial$ which is not natural
either. Below we will manage, however, to work out the procedure
sketched above within the framework of homomorphisms between
semi--holonomic jet modules. Thus the resulting operators $L$ will be
even strongly invariant. 

\subsection{}\label{5.2}
Each $P$--module ${\Bbb V}$ enjoys a decomposition
$$
\Bbb V=\Bbb V_{i_0}\oplus\Bbb V_{i_0+1}\oplus\dots\oplus\Bbb V_{i_0+k}
$$
as a $G_0$--module, where the submodules $\Bbb V_i$ are distinguished 
by the requirement that the grading element 
$E\in {\frak g}_0$ (cf. \ref{2.1}) 
 acts by scalar multiplication by $i$.  
The action of the elements $Z\in {\frak g_j}$
then maps $\Bbb V_i$ into $\Bbb V_{i+j}$ and so for each $j=0,\dots,k$
the subspace $\Bbb V^j:=\Bbb V_{i_0+j}\oplus\dots\oplus\Bbb V_{i_0+k}$
is a $P$--submodule of $\Bbb V$. In particular, this decomposition
of an irreducible $G$--module $\Bbb V$, viewed as $P$-module, runs from
$\Bbb V_{-k}$ to $\Bbb V_k$, where $\Bbb V_k$ is the $P$--submodule
generated by the highest weight of $\Bbb V$.
   
Now, let $\Bbb E_{i_0}$ be an irreducible component of $H^n(\frak
g_-,\Bbb V)$, on which the grading element acts by multiplication by
$i_0$. Then we can view $\Bbb E_{i_0}$ as a $G_0$ submodule of 
$\ker(\square)\subset\La^n\frak p_+\otimes\Bbb V$ and we write $\Bbb E$ for
the $P$--submodule in $\La^n\frak p_+\otimes\Bbb V$ generated by $\Bbb
E_{i_0}$.  Let
$$
\Bbb E=\Bbb E_{i_0}\oplus\dots\oplus\Bbb E_{i_0+r}
$$ 
be the above $G_0$--module decomposition 
according to eigenvalues of the grading element. Then 
the action of $\frak g_\ell$ maps each $\Bbb E_{i_0+i}$ to $\Bbb
E_{i_0+i+\ell}$. For each $i=1,\dots,r+1$ we have the $P$--submodule
$\Bbb E^i$ as above, so we can form the quotient $\Bbb E/\Bbb E^i$,
which is as a $G_0$--module isomorphic to $\Bbb
E_{i_0}\oplus\dots\oplus\Bbb E_{i_0+i-1}$. In particular, $\Bbb E/\Bbb
E^1$ is again the irreducible module $\Bbb E_{i_0}$ we started with
but now viewed as a $P$--module, and $\Bbb E/\Bbb E^{r+1}=\Bbb E$. 

\begin{lem*}
(1) $\Bbb E\subset\ker(\partial^*)$ and $\Bbb
E^1\subset\im(\partial^*)$.\newline 
(2) The Laplacian $\square$ restricts to a $G_0$--isomorphism $\Bbb
E_{i_0+i}\to\Bbb E_{i_0+i}$ for each $i=1,\dots,r$. 
\end{lem*}
\begin{proof}
(1) The first part is clear, since $\ker(\partial^*)$ is a $P$--submodule
which by construction contains $\Bbb E_{i_0}$. Since we have already
seen in Lemma \ref{4.4} 
that the action of $\frak p_+$ maps
$\ker(\partial^*)$ to $\im(\partial^*)$, the second part is also
clear.\newline 
(2) We have already noted in \ref{5.1} 
above that $\square$ restricts
to an automorphism on $\im(\partial^*)$. Hence it suffices to prove
that $\square(\Bbb E_{i_0+i})\subset\Bbb E_{i_0+i}$. By Corollary
\ref{3.9}, we have for all $e\in \Bbb E$, $Z\in {\frak g}_1$ 
$$
\partial(Z\cdot e)=Z\cdot\partial(e)-
(n+1)\sum_{|\eta_\al|=1}\eta_\al\wedge [Z,\xi_\al]\cdot e. 
$$ 
Applying $\partial^*$ to the first term we get
$Z\cdot\square(e)$. 

Let us first take
$e_0\in\Bbb E_{i_0}$, and consider $\square(Z\cdot 
e_0)=\partial^*(\partial(Z\cdot e_0))$. Then the first term vanishes while
each summand in the second term is contained
in $\partial^*(\frak g_1\wedge \frak g_0\cdot\Bbb
E_{i_0})\subset\partial^*(\frak
g_1\wedge\Bbb E_{i_0})$. Since $\Bbb E_{i_0}\subset\ker(\partial^*)$,
Lemma \ref{4.4} 
implies that $\partial^*(\frak g_1\wedge\Bbb
E_{i_0})\subset \frak g_1\cdot\Bbb E_{i_0}\subset\Bbb
E_{i_0+1}$. Thus, we see that $\square(\Bbb E_{i_0+1})\subset\Bbb
E_{i_0+1}$. Now one can proceed inductively in the same way to show
that $\square(\Bbb E_{i_0+i})\subset\Bbb E_{i_0+i}$. 
\end{proof}

\subsection{}\label{5.3}
The actual construction of the splitting operators is a little
tricky. The problem is that the individual steps in the construction
sketched in \ref{5.1} 
are induced by maps on jet--modules which are not $P$--module
homomorphisms but only restrict to $P$--module homomorphisms on
appropriate submodules, which also have to be constructed during the
procedure. 

For $j\geq i\ge0$ we denote by $\pi^j_i$ the canonical
projection $\Bbb E/\Bbb E^j\to\Bbb E/\Bbb E^i$, which is a
homomorphism of $P$--modules. Clearly, $\pi^i_i$ is the identity and
$\pi^j_i\o\pi^k_j=\pi^k_i$ for $i\leq j\leq k$. By $p_i:\Cal J^1(\Bbb
E/\Bbb E^i)\to\Bbb E/\Bbb E^i$ we denote the footpoint projection,
which is a $P$--homomorphism, too. For any element $\ps$ in a general 
$G_0$--module,  we denote by $\ps_i$ the component of $\ps$ on
which the grading element $E$ acts by multiplication by $i_0+i$. Note
that the mapping $\ps\mapsto \ps_i$ is only a $G_0$--homomorphism and
not a $P$--homomorphism. Finally, let us denote by $j_i:\Bbb E/\Bbb
E^i\to \Bbb E/\Bbb E^{i+1}$ the $G_0$--homomorphism given by the
inclusion $\Bbb E_{i_0}\oplus\dots\oplus\Bbb E_{i_0+i-1}\to\Bbb
E_{i_0}\oplus\dots\oplus\Bbb E_{i_0+i}$. Again, this is obviously {\em 
not\/} a $P$--homomorphism. Finally, let $\Alt:\frak
p_+\otimes\La^n\frak p_+\otimes\Bbb V\to\La^{n+1}\frak p_+\otimes\Bbb
V$ denote the alternation mapping. This is a $P$--homomorphism
preserving homogeneous degrees.

For $i=1,\dots,r+1$ consider now the module $\Cal J^1(\Bbb E/\Bbb
E^i)$. A typical element of this module is a pair $(e,\ps)$, with
$e\in\Bbb E/\Bbb E^i$ and 
$$
\ps\in\frak p_+\otimes\Bbb E/\Bbb E^i\subset
\frak p_+\otimes\La^n\frak p_+\otimes\Bbb V.
$$ 
Now we define a mapping $L_i:\Cal J^1(\Bbb E/\Bbb E^i)\to\Bbb E/\Bbb
E^{i+1}$ by 
$$
L_i(e,\ps)=j_i(e)-(n+1)\square^{-1}\partial^*((\Alt(\ps))_i).
$$
In particular, if $\ps=Z\otimes f$ for $Z\in\frak p_+$ and $f\in\Bbb
E/\Bbb E^i$, then $L_i(e,Z\otimes
f)=j_i(e)-(n+1)\square^{-1}\partial^*((Z\wedge f)_i)$. Now the main
technical step in the construction is the following

\begin{prop}\label{5.4}
The maps $L_i:\Cal J^1(\Bbb E/\Bbb E^i)\to\Bbb E/\Bbb E^{i+1}$ have
the following properties:\newline
(1) $L_i$ is a $G_0$--homomorphism and $\pi_i^{i+1}\o L_i=p_i$.\newline
(2) For $\Ps\in\Cal J^1(\Bbb E/\Bbb E^i)$ and $W\in\frak g_1$, we have
$$
L_i(W\cdot\Ps)-W\cdot L_i(\Ps)=
\square^{-1}\big(W\cdot
(\square\o j_i\o(L_{i-1}\o\Cal J^1(\pi_{i-1}^{i})-p_i)(\Ps))\big).
$$
In particular, $L_1$ is a $P$--homomorphism.
\end{prop}
\begin{proof}
(1) The fact that $L_i$ is a $G_0$--homomorphism follows immediately
    from the fact that $\Cal J^1(\Bbb E/\Bbb E^i)\cong\Bbb E/\Bbb
    E^i\oplus(\frak p_+\otimes\Bbb E/\Bbb E^i)$ as a $G_0$--module,
    see \ref{3.2} 
    and the definition of $L_i$. Moreover, since $\Alt$,
    $\partial^*$,
 and $\square$ all preserve homogeneities, the last
    term in the definition of $L_i$ is homogeneous of degree $i_0+i$,
    so it lies in the kernel of $\pi_i^{i+1}$ and the second part
    follows.

(2) Clearly, it suffices to check this for elements $\Ps$ of the form
    $(e,Z\otimes f)$ with $e,f\in\Bbb E/\Bbb E^i$ and $Z\in\frak
    p_+$. By definition of the action on jets, see \ref{3.2}, 
    we see that $W\cdot (e,Z\otimes f)$ has footpoint $W\cdot e$, while the
    homogeneous part of degree $i_0+i$ of the second component is
    given by 
$$
\sum_{|\eta_\al|=1}\eta_\al\otimes [W,\xi_\al]\cdot
    e_{i-1}+W\cdot(Z\otimes f)_{i-1}.
$$
Consequently, 
\begin{align*}
L_i(W\cdot(e,Z\otimes f))=\ &j_i(W\cdot e)-(n+1)\square^{-1}\partial^*(
\sum_{|\eta_\al|=1}\eta_\al\wedge [W,\xi_\al]\cdot e_{i-1})-\\
&-(n+1)\square^{-1}\partial^*(W\cdot(Z\wedge f)_{i-1}).
\end{align*}
By Corollary \ref{3.9} 
the second term on the right hand side of this equation just gives 
$$
\square^{-1}\partial^*(\partial(W\cdot e_{i-1})-
W\cdot\partial(e_{i-1}))=W\cdot
e_{i-1}-\square^{-1}(W\cdot\square(e_{i-1})),
$$
where we have used that $\partial^*$ is a $P$--homomorphism, $e_{i-1}$
and $W\cdot e_{i-1}$ lie in the kernel of $\partial^*$, and that we
are in a region where the Laplacian is invertible. On the other hand,
we clearly have $j_i(W\cdot e)+W\cdot e_{i-1}=W\cdot j_i(e)$, since
$W\in\frak g_1$ and $e_{i-1}$ is the highest nonzero homogeneous
component of $e$. Finally, we clearly have $W\cdot L_i(e,Z\otimes
f)=W\cdot j_i(e)$, since the rest lies in the component of maximal
homogeneity, on which $\frak p_+$ acts trivially. Thus, we have
arrived at
\begin{multline*}
L_i(W\cdot(e,Z\otimes f))-W\cdot L_i(e,Z\otimes f)=\\
=-\square^{-1}(W\cdot\square(e_{i-1}))-(n+1)
\square^{-1}(W\cdot\partial^*((Z\wedge f)_{i-1})),
\end{multline*}
where we have used once more the fact that $\partial^*$ is a
$P$--homomorphism. 

On the other hand, consider $\Cal J^1(\pi_{i-1}^i)(e,Z\otimes f)$. The
footpoint of this element is just $\pi_{i-1}^i(e)$, while in the jet
part, the component of maximal homogeneity must coincide with
$(Z\otimes f)_{i-1}$. Consequently, we get 
$$
L_{i-1}(\Cal J^1(\pi_{i-1}^i)(e,Z\otimes
f))=j_{i-1}(\pi_{i-1}^i(e))-(n+1)\square^{-1}\partial^*((Z\wedge
f)_{i-1}). 
$$
Subtracting $e=p_i(e,Z\otimes f)$ from this, we get
$$
-e_{i-1}-(n+1)\square^{-1}\partial^*((Z\wedge f)_{i-1}),
$$
and the formula follows. In the case $i=1$, we get
$L_1(W\cdot\Ps)-W\cdot L_1(\Ps)=-\square^{-1}(W\cdot(\square\o j_1\o
p_1)(\Ps))$, which vanishes, since $\Bbb
E_{i_0}\subset\Ker(\square)$. Hence, $L_1$ is equivariant for the
action of $\frak g_1$ and thus a $P$--homomorphism by (1) and \ref{2.4}. 
\end{proof}

\subsection{}\label{5.5}
Now we inductively define subspaces $\tilde\Cal J^1(\Bbb E/\Bbb
E^i)\subset \Cal J^1(\Bbb E/\Bbb E^i)$ by $\tilde\Cal J^1(\Bbb E/\Bbb 
E^1)=\Cal J^1(\Bbb E/\Bbb E^1)$ and 
$$
\tilde\Cal J^1(\Bbb E/\Bbb E^{i+1})=\Cal J^1(\pi_i^{i+1})^{-1}(\tilde\Cal
J^1(\Bbb E/\Bbb E^i))\cap \Ker(L_i\o\Cal J^1(\pi_i^{i+1})-p_{i+1}). 
$$

\begin{prop*}
For each $i=1,\dots,r+1$ the space $\tilde\Cal J^1(\Bbb E/\Bbb
E^i)\subset\Cal J^1(\Bbb E/\Bbb E^i)$ is a $P$--submodule and $L_i$
restricts to a homomorphism $\tilde\Cal 
J^1(\Bbb E/\Bbb E^i)\to\Bbb E/\Bbb E^{i+1}$ of $P$--modules. 
Moreover, for each $k<i$ we have 
$$
\Cal J^1(\pi^i_k)\big(\tilde\Cal J^1(\Bbb E/\Bbb
E^i)\big)\subset\tilde\Cal J^1(\Bbb E/\Bbb E^k),
$$
and on $\tilde\Cal J^1(\Bbb E/\Bbb E^i)$ we have $\pi^i_{k+1}\o
p_i=L_k\o\Cal J^1(\pi^i_k)$.  
\end{prop*}
\begin{proof}
For $i=1$ the first two properties are satisfied by definition of
$\tilde\Cal J^1(\Bbb E/\Bbb E^1)$ and Proposition \ref{5.4}(2), while
the last two properties are trivially satisfied. If we inductively
assume that the result has been proved for $i-1$, then $\Cal
J^1(\pi^i_{i-1})^{-1}(\tilde\Cal J^1(\Bbb E/\Bbb E^{i-1}))$ is a
$P$--submodule of $\Cal J^1(\Bbb E/\Bbb E^i)$, and $L_{i-1}\o\Cal
J^1(\pi^i_{i-1})-p_i$ defines a $P$--homomorphism from this submodule
to $\Bbb E/\Bbb E^i$, so $\tilde\Cal J^1(\Bbb E/\Bbb E^i)$ is a
$P$--submodule. Moreover, Proposition \ref{5.4}(2) immediately implies
that the restriction of $L_i$ to this submodule is equivariant under
the action of $\frak g_1$ and thus $L_i$ restricts to a
$P$--homomorphism on that submodule by Proposition \ref{5.4}(1) and
\ref{2.4}. Moreover, we get the last two properties in the case
$k=i-1$. 

For $k<i-1$, note first that
$\pi^i_k=\pi^{i-1}_k\o\pi^i_{i-1}$ immediately implies that $\Cal
J^1(\pi^i_k)\big(\tilde\Cal J^1(\Bbb E/\Bbb E^i)\big)\subset\tilde\Cal
J^1(\Bbb E/\Bbb E^k)$ by induction. Finally, we compute
\begin{align*}
L_k\o\Cal J^1(\pi^i_k)&=L_k\o\Cal J^1(\pi^{i-1}_k)\o\Cal
J^1(\pi^i_{i-1})=\pi^{i-1}_{k+1}\o p_{i-1}\o\Cal J^1(\pi^i_{i-1})=\\
&=\pi^{i-1}_{k+1}\o\pi^i_{i-1}\o p_i=\pi^i_{k+1}\o p_i,
\end{align*}
by functoriality of $\Cal J^1$, induction, and the definition of the
jet prolongation of a homomorphism. 
\end{proof}

For $k\geq 2$ and $i=1,\dots, r+1$ we inductively define 
$$
\tilde\Cal J^k(\Bbb E/\Bbb E^i):=\Cal J^1(\tilde\Cal J^{k-1}(\Bbb
E/\Bbb E^i))\cap\bar\Cal J^k(\Bbb E/\Bbb E^i).
$$
By Proposition \ref{5.5} and \ref{3.4} 
it follows inductively that $\tilde\Cal J^k(\Bbb E/\Bbb E^i)$ is a
$P$--submodule in both modules on the right hand side of the
definition. For $i=1$, we obtain $\tilde\Cal J^k(\Bbb
E/\Bbb E^1)=\bar\Cal J^k(\Bbb E/\Bbb E^1)$, so we simply get the full
semiholonomic jet--module in this case. Moreover, a simple inductive
argument shows for all $\ell<k$, and $i$ 
$$
\tilde \Cal J^k(\Bbb E/\Bbb E^i) \subset \bar \Cal J^\ell(\tilde\Cal
J^{k-\ell}(\Bbb E/\Bbb E^i))\cap \bar\Cal J^k(\Bbb E/\Bbb E^i)
.$$ 

For each of the homomorphisms $L_i:\tilde\Cal J^1(\Bbb E/\Bbb
E^i)\to\Bbb E/\Bbb E^{i+1}$ we can now restrict the semiholonomic jet 
prolongation $\bar\Cal J^k(L_i)$ to the submodule $\tilde\Cal
J^{k+1}(\Bbb E/\Bbb E^i)\subset \bar \Cal J^{k}(\tilde\Cal J^1(\Bbb E/\Bbb E^i))$ to
obtain a $P$--homomorphism 
$$
\bar\Cal J^k(L_i):\tilde\Cal J^{k+1}(\Bbb
E/\Bbb E^i)\to \bar\Cal J^k(\Bbb E/\Bbb E^{i+1})
.$$

\begin{thm}\label{5.6} Let ${\Bbb E}_{i_0}$ be an irreducible component in
the cohomology $\Bbb H^n_{\Bbb V}$ which generates the  $P$--submodule 
${\Bbb E}$ in $\La^n{\frak p}_+\otimes {\Bbb V}$, cf. \ref{5.1}.
For each $k\geq 1$ and $i=1,\dots,r+1$ we have
$$
\bar\Cal J^k(L_i)\big(\tilde\Cal J^{k+1}(\Bbb E/\Bbb E^i)\big)\subset
\tilde\Cal J^k(\Bbb E/\Bbb E^{i+1}).
$$
In particular, the composition 
$$
L:=L_r\o\bar\Cal J^1(L_{r-1})\o\dots\o\bar\Cal J^{r-1}(L_1)
$$ defines
a $P$--homomorphism $L:\bar\Cal J^r(\Bbb E/\Bbb E_1)\to\Bbb E$. Since
by definition $\Bbb E$ is a $P$--submodule of $\La^n\frak
p_+\otimes\Bbb V$, this homomorphism induces a strongly invariant
operator $\Ga(E_{i_0}M)\to\Ker(\partial^*)\subset\Om^n(M;VM)$, which
splits the algebraic projection $\Ker(\partial^*)\to\Ga(E_{i_0}M)$
described in \ref{5.1}. 
\end{thm}
\begin{proof}
Let us first consider the case $k=1$. So we have to show that $\Cal
J^1(L_i)\big(\tilde\Cal J^2(\Bbb E/\Bbb E^i)\big)\subset\tilde\Cal
J^1(\Bbb E/\Bbb E^{i+1})$. By definition of $\tilde\Cal J^1(\Bbb
E/\Bbb E^{i+1})$, we first have to consider the composition $\Cal
J^1(\pi^{i+1}_i)\o\Cal J^1(L_i)=\Cal J^1(\pi^{i+1}_i\o L_i)$. By Proposition
\ref{5.4}(1), this equals $\Cal J^1(p_i)$. Since $\tilde\Cal J^2(\Bbb
E/\Bbb E^i)\subset\bar\Cal J^2(\Bbb E/\Bbb E^i)$, this projection
coincides with the restriction of the canonical projection $\bar\Cal
J^2(\Bbb E/\Bbb E^i)\to\Cal J^1(\Bbb E/\Bbb E^i)$, and since
$\tilde\Cal J^2(\Bbb E/\Bbb E^i)\subset\Cal J^1(\tilde\Cal J^1(\Bbb
E/\Bbb E^i))$, this canonical projection has values in $\tilde\Cal
J^1(\Bbb E/\Bbb E^i)$. Thus, we have verified that 
$\Cal J^1(L_i)\big(\tilde\Cal J^2(\Bbb E/\Bbb E^i)\big)\subset\Cal
J^1(\pi^{i+1}_i)^{-1}\big(\tilde\Cal J^1(\Bbb E/\Bbb E^i)\big)$.

But then it also follows that $L_i\o\Cal J^1(\pi^{i+1}_i)\o\Cal J^1(L_i)$
coincides with the composition of $L_i$ with the canonical projection
$\tilde\Cal J^2(\Bbb E/\Bbb E^i)\to\tilde\Cal J^1(\Bbb E/\Bbb E^i)$,
which by definition of the jet prolongation of a homomorphism 
(see \ref{3.2}) 
coincides with $p_{i+1}\o \Cal J^1(L_i)$ and the proof in the
case $k=1$ is complete.

The case $k\geq 2$ now immediately follows from the definitions by
induction. Thus, also the existence of $L$ and the corresponding
strongly invariant operator is clear. The fact that this operator
splits the algebraic projection follows from the fact that by Lemma
\ref{5.2}(1) this algebraic projection is induced by the canonical
projection $\Bbb E\to\Bbb E/\Bbb E^1$ and the fact that $\pi_i^{i+1}\o
L_i=p_i$ from Proposition \ref{5.4}(1). 
\end{proof}

Next, we consider the composition of $d_\Bbb V$ with the
operator corresponding to $L$. The corresponding homomorphism on jet
modules can be computed as the restriction to $\bar\Cal J^{r+1}(\Bbb
E/\Bbb E^1)$ of $d_\Bbb V\o\Cal J^1(L)$. 

\begin{prop}\label{5.5a} For each irreducible $G$-module $\Bbb V$, and
irreducible $G_0$--submodule $\Bbb E_{i_0}\subset \Bbb H^n_\Bbb V=
H^n({\frak g}_-,\Bbb V)$, the composition 
$$
d_{\Bbb V}\o \Cal J^1 
L:\bar\Cal J^{r+1}(\Bbb E/\Bbb E^1)\to \La^{n+1}\frak p_+\otimes \Bbb V
$$
has values in $\ker\partial ^*$. The composition with the 
projection to the cohomology $\pi_H:(\La^n{\frak p}_+\otimes \Bbb V)\cap
(\ker\partial^*)\to\Bbb H^{n+1}_\Bbb V=H^{n+1}({\frak g}_-,\Bbb V)$ gives the
$P$-module homomorphism 
$$
\pi_H\o d_{\Bbb V}\o \Cal J^1L:\bar\Cal J^{r+1}(\Bbb E/\Bbb E^1)\to
\Bbb H^{n+1}_\Bbb V. 
$$ 
For each $n=0,\dots,\operatorname{dim} M-1$, there is the strongly invariant
differential operator 
$$
D^{\Bbb V}:\Gamma(H^n_{\Bbb V}M)\to\Gamma(H^{n+1}_{\Bbb V}M)
$$ 
whose restrictions to the subbundles $E_0M$ are determined by the above
$P$--module homomorphisms $\bar \Cal J^{r+1}(\Bbb E/\Bbb E^1)\to \Bbb
H^{n+1}_{\Bbb V}$. 
\end{prop}
\begin{proof}
Consider first the map $\partial^*\o d_\Bbb V:\Cal J^1(\Bbb
E)\to\La^n\frak p_+\otimes\Bbb V$. By definition of $d_\Bbb V$, 
Lemma \ref{4.4}, 
and using the fact that $\Bbb E\subset\Ker(\partial^*)$, we see that
this maps $(e,Z\otimes f)\in\Cal J^1(\Bbb E)$ to
$\partial^*\partial(e)+(n+1)\partial^*(Z\wedge
f)=\square(e)-(n+1)Z\cdot f$, so $\partial^*\o d_\Bbb V:\Cal J^1(\Bbb 
E)\to\Bbb E$. Now Theorem \ref{5.6} applied to $\Cal J^1(L_r)$
shows, that $\Cal J^1(L)$ has values in the submodule $\tilde\Cal
J^1(\Bbb E)\subset\Cal J^1(\Bbb E)$, and we claim that $\partial^*\o
d_\Bbb V$ restricts to zero on that submodule. 

To simplify notations, let us write $p:\Cal J^1(\Bbb E)\to\Bbb E$ for
the footpoint projection $p_{r+1}$ and $\pi_i$ for $\pi^{r+1}_i$. For
$i\leq r+1$ consider the $P$--homomorphism $\pi_i\o\partial^*\o d_\Bbb
V:\Cal J^1(\Bbb E)\to\Bbb E/\Bbb E^i$. By definition, this maps
$(e,Z\otimes f)$ to $\pi_i(\square(e))+(n+1)\pi_i(\partial^*(Z\wedge
f))$. Since the Laplacian and $\partial^*$ both preserve homogeneous
degrees, we may rewrite the first summand as $\square(\pi_i(e))$ and
the second summand as
$(n+1)\pi_i(\partial^*(Z\wedge\pi_{i-1}(f)))$. 

On the other hand, consider $\Cal J^1(\pi_{i-1}):\Cal J^1(\Bbb
E)\to\Cal J^1(\Bbb E/\Bbb E^{i-1})$. This maps $(e,Z\otimes f)$ to
$(\pi_{i-1}(e),Z\otimes\pi_{i-1}(f))$, and applying $L_{i-1}$ to this
element, we get
$j_{i-1}(\pi_{i-1}(e))-(n+1)\square^{-1}\partial^*((Z\wedge\pi_{i-1}(f))_i)$. 
Finally, $\pi_i\o p$ maps $(e,Z\otimes f)$ to
$\pi_i(e)$. Consequently, $\square\o(\pi_i\o p-L_{i-1}\o\Cal
J^1(\pi_{i-1}))$ maps $(e,Z\otimes f)$ to 
$$
\square(\pi_i(e))-j_{i-1}(\square(\pi_{i-1}(e)))+
(n+1)\partial^*((Z\wedge\pi_{i-1}(f))_i),
$$
and the last summand in this expression equals
$$
(\pi_i-j_{i-1}\o\pi_{i-1})((n+1)\partial^*(Z\wedge f)),
$$ 
since $\partial^*$ preserves homogeneous degrees. Hence, we see that
on $\Cal J^1(\Bbb E)$ we get the equation
$$
\pi_i\o\partial^*\o d_\Bbb V-j_{i-1}\o\pi_{i-1}\o\partial^*\o d_\Bbb
V= \square\o\big(\pi_i\o p-L_{i-1}\o\Cal J^1(\pi_{i-1})\big).
$$
In fact, this equation is exactly what we were aiming at in the
motivation for the whole construction in \ref{5.1}. 
But on the submodule $\tilde\Cal J^1(\Bbb E)$, the right hand side of
the above formula vanishes identically by Proposition \ref{5.5}. Thus,
iterated application of this formula shows that on $\tilde\Cal
J^1(\Bbb E)$ we have
$$
\partial^*\o d_\Bbb V=\pi_{r+1}\o\partial^*\o d_\Bbb
V=j_r\o\pi_r\o\partial^*\o d_\Bbb V=\dots=j_1\o\pi_1\o \partial^*\o d_\Bbb V.
$$
But $\pi_1\o\partial^*\o d_\Bbb V$ maps $(e,Z\otimes f)$ to
$\square(\pi_1(e))$, which vanishes since $\Bbb E_{i_0}$ is contained
in the kernel of the Laplacian, so we have proved $\partial^*\o d_\Bbb
V\o\Cal J^1(L)=0$. All the rest is now an immediate consequence. 
\end{proof}

\subsection{Definition}\label{5.6a} 
Let $(\Cal G,\om)$ be a (real or complex) parabolic
geometry on a manifold $M$.
The construction above has given rise to a sequence of strongly
invariant operators $D^{\Bbb V}$
$$
\xymatrix@1{
{0}\ar[r] & {\Ga(H^0_\Bbb VM)} \ar[r]^-{D^{\Bbb V}} & {\Ga(H^1_\Bbb VM)} 
\ar[r]^-{D^{\Bbb V}}&\dots
\ar[r]^-{D^{\Bbb V}}&{\Ga(H^{\text{dim}(G/P)}_\Bbb VM)}\ar[r] &{0}}.
$$ 
which is called the \emph{Bernstein--Gelfand--Gelfand sequence\/} or
\emph{BGG--se\-quence\/} determined by the $G$-module $\Bbb V$. 
 
All bundles in this sequence correspond to completely reducible
representations of $P$, so they all split into  direct sums of bundles
corresponding to irreducible representations. Let us also remark that the
construction applies to both real and complex settings. 
Next, we will show that
in the flat case this sequence is a resolution of the constant sheaf
$\Bbb V$. Since by Kostant's version of the Bott--Borel--Weil theorem,
the bundles occurring in this resolution in the complex case are
exactly the same bundles as in the Bernstein--Gelfand--Gelfand
resolution, we have obtained curved analogs of this resolution even in
the real case. 

The main step towards the proof that we often get a resolution is formulated
in the next lemma for the general real curved case. For the complex analog
see below.

\begin{lem}\label{5.7} Let $(\Cal G,\om)$ be a real parabolic geometry on a
manifold $M$ and 
let $s\in\Om^n(M;VM)$ be a $VM$--valued $n$--form. Then:\newline
(1) There is an element $t\in\Om^{n-1}(M;VM)$ such that $s+d_\Bbb
V(t)$ lies in $\ker(\partial^*)$.\newline 
(2) If $s$ and $d_\Bbb V(s)$ both lie in $\ker(\partial^*)$, then
$s=L(\pi_H(s))$.\newline
(3) If $d_\Bbb V^2(\ker(\partial^*))\subset\ker(\partial^*)$, then the
diagram  
\begin{equation*}
\xymatrix@R=7pt{
& {\Om^0(M;VM)} \ar[r]^-{d_{\Bbb V}}
& {\Om^1(M;VM)} \ar[r]^-{d_{\Bbb V}}
& {\dots}
\\
{\Bbb V} \ar[ur] \ar[dr] 
\\
& {\Gamma(H^0_{\Bbb V}M)} \ar[uu]_L \ar[r]^-{D^{\Bbb V}} 
& {\Gamma(H^1_{\Bbb V}M)} \ar[uu]_L \ar[r]^-{D^{\Bbb V}}
&{\dots}
}
\end{equation*}
is commutative. In particular, $D^{\Bbb V}\o D^{\Bbb V}=0$ whenever
$d_{\Bbb V}\o d_{\Bbb V}=0$.
\end{lem}
\begin{proof}
(1) Put $\Cal G_0=\Cal G/P_+$ and choose a global $G_0$--equivariant
section $\si:\Cal G_0\to\Cal G$ as indicated in \ref{2.14}. 
Then we get a smooth
map $\tau:\Cal G\to P_+$ characterized by $u=\si(p(u))\cdot\tau(u)$
for all $u\in\Cal G$, and $u\mapsto (p(u),\tau(u))$ is a
diffeomorphism $\Cal G\to \Cal G_0\x P_+$. Using this, we get an
isomorphism (depending on $\si$) between $\Om^n(M;VM)$ and the space
of smooth $G_0$--equivariant functions $\Cal G_0\to\La^n\frak
p_+\otimes\Bbb V$. But $\square^{-1}\o\partial^*$ is a
$G_0$--homomorphism $\La^n\frak p_+\otimes\Bbb V\to\La^{n-1}\frak
p_+\otimes\Bbb V$ such that
$e-\partial(\square^{-1}\o\partial^*(e))\in\ker(\partial^*)$ for all 
$e\in \La^n{\frak p_+\otimes \Bbb V}$. 

Now, let $\underline s:\Cal G_0\to \La^n\frak p_+\otimes\Bbb V$ be the
$G_0$--equivariant map corresponding to the lowest homogeneous
component $s_j$ of the given $n$--form $s$
such that $\partial^*(s_j)\neq 0$. Passing from
$-\square^{-1}\o\partial^*\o \underline s$ back to a $P$--equivariant
map $t:\Cal G\to \La^{n-1}\frak p_+\otimes\Bbb V$, we see that the
homogeneous components up to degree $j$ of $\partial^*(s+d_\Bbb V(t))$
vanish on the image of $\si$ and thus on the whole $\Cal G$ by
equivariancy. Inductively, we can find an element 
$t$ with the required properties.

\noindent
(2) Put $s'=\pi_H(s)$. By construction of the operators $L$, we know
that $L(s')\in\ker(\partial^*)$, $\pi_H(L(s'))=s'$, and $d_\Bbb
V(L(s'))\in\ker(\partial^*)$. Thus, we see that
$s-L(s')\in\im(\partial^*)$ and $d_\Bbb
V(s-L(s'))\in\ker(\partial^*)$. Let $a_j$ be the lowest possibly nonzero
homogeneous component of $s-L(s')$. Then the lowest possibly nonzero
component of $d_\Bbb V(s-L(s'))$ is $\partial(a_j)$. Since
$\ker(\partial^*)$ is complementary to $\im(\partial)$ we must have
$\partial(a_j)=0$. But on the other hand, we know that
$a_j\in\im(\partial^*)$ which is complementary to $\ker(\partial)$, so
we must have $a_j=0$.

\noindent(3) For $s\in\Ga(H^n_\Bbb VM)$, consider the element $d_\Bbb
V(L(s))\in\Om^{n+1}(M;VM)$. By Proposition \ref{5.5a}, this lies in
$\ker(\partial^*)$. Moreover, since $L(s)\in\ker(\partial^*)$, our
assumption on $d_\Bbb V^2$ implies that $d_\Bbb V(d_\Bbb
V(L(s)))\in\ker(\partial^*)$. Hence from (2) we get 
$d_\Bbb V(L(s))=L(\pi_H(d_\Bbb V(L(s))))=L(D^\Bbb V(s))$. 

The last claim is obvious.
\end{proof}

\begin{lem}\label{5.7a} Let $(\Cal G,\om)$ be a complex parabolic geometry
on a complex manifold $M$. Then the second and third assertion in Lemma 
\ref{5.7} remain valid with the same assumptions, while the claim
\ref{5.7}(1) holds true under the additional assumption that the holomorphic
bundle  $\Cal G\to\Cal G_0$ admits a global holomorphic $G_0$-equivariant
section. This additional requirement is always fulfilled locally.
\end{lem}
\begin{proof}
All arguments in the proof of (2) and (3) in \ref{5.7} are on the level of
the $P$-modules and so they go equally through for both real and complex
settings. The only difference in (1) is the argument which constructs the
global section by means of the smooth partition of unity. Once we assume the
existence of the global section, the rest is clear again. Now, any
point in $M$ has an open neighborhood $U\subset M$ such that both $\Cal
G$ and $\Cal G_0$ are trivial over $U$. Since $G_0\x P_+$ and $P$ are
diffeomorphic, and the map in one direction is obviously holomorphic,
they are biholomorphic. Thus, the complex parabolic geometry $\Cal
G|_U\to U$ admits appropriate global holomorphic $G_0$--equivariant
section.
\end{proof}

\begin{thm}\label{5.8}
Let $(\Cal G,\om)$ be a real parabolic geometry of the type $(G,P)$ 
on a manifold $M$, $\Bbb V$ be a $G$--module.  If the twisted de~Rham
sequence
$$
\xymatrix@1@C=5mm{
{0}\ar[r] & {\Om^0(M;VM)} \ar[r]^-{d_{\Bbb V}} & {\Om^1(M; VM)} 
\ar[r]^-{d_{\Bbb V}}&\dots
\ar[r]^-{d_{\Bbb V}}&{\Om^{\operatorname{dim}(G/P)}(M;VM)}\ar[r] &{0}}.
$$
is a complex, then also the Bernstein--Gelfand--Gelfand
sequence 
$$
\xymatrix@1{
{0}\ar[r] & {\Ga(H^0_\Bbb VM)} \ar[r]^-{D^{\Bbb V}} & {\Ga(H^1_\Bbb VM)} 
\ar[r]^-{D^{\Bbb V}}&\dots
\ar[r]^-{D^{\Bbb V}}&{\Ga(H^{\operatorname{dim}(G/P)}_\Bbb VM)}\ar[r] &{0}}
$$
defined in \ref{5.6a} is a complex, and they both compute the same
cohomology.

The same statement is true for complex parabolic geometries $(\Cal G,\om)$ 
under the additional requirement that $\Cal G\to \Cal G_0=\Cal G/P_+$ admits
a global holomorphic $G_0$--equivariant section. 
\end{thm}

\subsection*{Remark}
In particular, the complex version of the Theorem may be reformulated
as follows: {\em If the twisted de~Rham sequence induces a complex on
the sheaf level, then the same is true for the
Bernstein--Gelfand--Gelfand sequence}. In particular, if the twisted
de~Rham sequence induces a resolution of $\Bbb V$, then so does the
BGG--sequence.

Now, the original
representation theoretical version of the (generalized) BGG--resolution
follows immediately by duality.
Moreover, let us notice that the global $G_0$--equivariant section as
required in the Theorem always exists over a dense open submanifold in
the homogeneous space $G/P$ (the so called big cell).
 
\begin{proof}
As we saw in Lemma \ref{5.7}, the BGG--sequence forms a complex whenever the
twisted de~Rham does. So let us assume, we deal with complexes. 

Since $d_\Bbb V^2=0$, \ref{5.7}(3) implies that $L$ is a morphism of the
corresponding complexes, hence the mapping 
$$\Gamma(H^n_{\Bbb V}M)\ni s' \mapsto L(s')\in \Om^n(M;VM)
$$
induces a morphism between the cohomologies.

Next, suppose that $s\in\Om^n(M;VM)$, $n\ge 1$ is such that $d_\Bbb
V(s)=0$. Then by \ref{5.7}(1) we find an element
$t\in\Om^{n-1}(M;VM)$ such that $s+d_\Bbb
V(t)\in\ker(\partial^*)$. But then $d_\Bbb V(s+d_\Bbb V(t))=0$ so by
\ref{5.7}(2) we know that $s+d_\Bbb V(t)=L(\pi_H(s+d_\Bbb
V(t)))$, and thus the mapping defined above is surjective. 

Finally, let us assume that $s'\in\Ga(H^{n-1}_\Bbb VM)$ is such that
there exists a $t\in\Om^{n-1}(M;VM)$ with $d_\Bbb V(t)=L(s')$. Then by
\ref{5.7}(1) we may without loss of generality assume that
$t\in\ker(\partial^*)$. But by assumption $d_\Bbb V(t)=L(s')$, so this
is also contained in $\ker(\partial^*)$, and hence $t=L(\pi_H(t))$ by
\ref{5.7}(2), and thus $L(s')=d_\Bbb V(L(\pi_H(t)))$ and
applying $\pi_H$ on both sides we get $s'=D^\Bbb V(\pi_H(t))$, and so
we get an isomorphism in the cohomology groups. 
\end{proof}

\begin{kor}\label{5.9}
Let $(\Cal G,\om)$ be a torsion free real parabolic geometry of type $(G,P)$ 
on $M$. Then the de~Rham cohomology of $M$ with coefficients in $\Bbb
K=\Bbb R$ or $\Bbb C$ is computed by the (much smaller) complex
$$
\xymatrix@1{
{0}\ar[r] & {\Ga(H^0_\Bbb KM)} \ar[r]^-{D^{\Bbb K}} & {\Ga(H^1_\Bbb KM)} 
\ar[r]^-{D^{\Bbb K}}&\dots
\ar[r]^-{D^{\Bbb K}}&{\Ga(H^{\operatorname{dim}(G/P)}_\Bbb KM)}\ar[r] &{0}}
.$$ 
\end{kor}
\begin{proof}
The covariant exterior differential corresponding to the choice of the
trivial $P$--module ${\Bbb K}$ coincides with the usual exterior
differential $d$. According to Lemma \ref{3.10}, the exterior covariant
differential coincides with our twisted exterior differential for all
torsion--free geometries. Thus the statement follows from \ref{5.8}.
\end{proof}

\begin{kor}\label{5.9a}
Let $(\Cal G,\om)$ be a flat real parabolic geometry. Then for any
representation $\Bbb V$ of $G$ the BGG--sequence 
$$
\xymatrix@1{
{0}\ar[r] & {\Ga(H^0_\Bbb VM)} \ar[r]^-{D^{\Bbb V}} & {\Ga(H^1_\Bbb VM)} 
\ar[r]^-{D^{\Bbb V}}&\dots
\ar[r]^-{D^{\Bbb V}}&{\Ga(H^{\operatorname{dim}(G/P)}_\Bbb VM)}\ar[r] &{0}}
$$
is a complex, which computes the twisted de~Rham cohomology of $M$
with coefficients in the bundle $VM$, which is defined as the
cohomology of the complex given by the covariant exterior derivative
with respect to the linear connection on $VM$ induced by the Cartan
connection $\om$, see \ref{3.10}.
\end{kor}
The importance of this corollary lies in the fact that while flat
parabolic geometries are locally isomorphic to the homogeneous model
$G/P$, they may be very different from $G/P$ from a global point of
view. Just keep in mind the broad variety of smooth manifolds
admitting a locally conformally flat Riemannian metric. In particular,
the bundle $VM$ is not trivial in general, so the twisted de~Rham
cohomology is a less trivial object than in the homogeneous case. 

On the other hand, we may always consider the obvious flat parabolic
geometry on the trivial $P$--bundle over $\Bbb
R^{\operatorname{dim}(G/P)}\cong\frak g_-$. In this case, the twisted
de~Rham cohomologies are obviously zero, so Corollary \ref{5.9a}
provides global resolutions of the constant sheaf $\Bbb V$ in this
case. Simple instances of such sequences are of basic importance in
various areas of mathematics, see for example \cite{MikeBGG}.  

\subsection{Remark}\label{5.10} 
As we have seen already, the $P$--modules ${\Bbb
H}^n_{\Bbb V}$ are completely reducible and so the natural 
bundles $H^n_{\Bbb V}M$ decompose into direct sums of irreducible bundles.
Consequently, also the operators $D^{\Bbb V}$ split into sums of 
operators between the irreducible natural bundles. In the case of the
homogeneous bundles, the latter operators (and sometimes also their
non-trivial compositions) are usually referred to as \emph{standard invariant
operators}. In particular, our construction provides a distinguished curved
analog for each of those standard operators. 

As we have underlined already 
in the introduction, no deep representation theoretical results had to be
applied in the construction of the BGG--sequences and in the proof of
Theorem \ref{5.8}. On the other hand, the full information of the Kostant's
version of Bott--Borel--Weil theorem on the Lie algebra cohomologies is
strictly necessary in order to get more explicit information about the
individual standard operators and the overall structure of the BGG--sequence
in the flat case. Moreover, further non-trivial operators with curvature
contributions in their symbols may appear in general. 

Let us also remark that the explicit formulae for the standard operators
were given in closed form in the terms of the underlying linear connections
on $M$ in \cite{CSS3} for all parabolic geometries with irreducible tangent
bundles, i.e. for all cases with $|1|$--graded Lie algebra ${\frak g}$. We
believe that the technique developed there should be extendible to the
general case, too.
 
\subsection{Remark}\label{5.11}
In the flat case, the twisted de~Rham complex can be viewed as a
filtered complex with the filtration given by homogeneous degrees. The
fact that the lowest homogeneous component of $d_\Bbb V$ is just
$\partial$ implies that the differential on the associated graded
complex is exactly $\partial$. Associated to this filtration there is
a spectral sequence which obviously converges and computes the twisted
de~Rham cohomology. Now from the construction of the operators
$D^\Bbb V$ it is obvious that when passing to the appropriate
subquotients, they induce the higher differentials in this spectral
sequence. Usually, these higher differentials are only well defined on
the appropriate subquotients, but due to the fact that we have a
(fairly simple) Hodge structure on the associated graded complex, we
can get a global definition in our setting.

\section{Example}\label{5}

\newcommand{\dynkin}[3]{%
\begin{picture}(50,15)
\put(5,2.5){\line(1,0){40}}
\put(5,2.5){\makebox(0,0){$\times$}}
\put(5,10){\hbox to0pt{\hss$\ssize #1$\hss}}
\put(25,2.5){\makebox(0,0){$\bullet$}}
\put(25,10){\hbox to0pt{\hss$\ssize #2$\hss}}
\put(45,2.5){\makebox(0,0){$\times$}}
\put(45,10){\hbox to0pt{\hss$\ssize #3$\hss}}
\end{picture}%
}

\newcommand{\dynkinG}[3]{%
\begin{picture}(50,15)
\put(5,2.5){\line(1,0){40}}
\put(5,2.5){\makebox(0,0){$\bullet$}}
\put(5,10){\hbox to0pt{\hss$\ssize #1$\hss}}
\put(25,2.5){\makebox(0,0){$\bullet$}}
\put(25,10){\hbox to0pt{\hss$\ssize #2$\hss}}
\put(45,2.5){\makebox(0,0){$\bullet$}}
\put(45,10){\hbox to0pt{\hss$\ssize #3$\hss}}
\end{picture}%
}

\newcommand{\SSS}[3]{\Cal S^{#1}_{[#2,#3]}}

We shall illustrate the power of our results in the simple case of 
5--dimensional partially integrable almost CR--structures, cf. Example 
\ref{2.10}. We believe that this simple geometry reflects many of the
general features of parabolic geometries and we can still write down the
whole BGG--sequences very explicitly at the same time. We hope that based on
this example, the reader is able to imagine the vast amount of invariant
operators which our main theorems produce for all parabolic
geometries. 

Let $M$ be a smooth manifold of odd dimension $2n+1$ together with a
distinguished rank $n$ complex subbundle $T^{CR}M$ of the tangent
bundle $TM$. Then the Lie bracket of vector fields induces a
skew--symmetric bundle map $\Cal L^\Bbb R:T^{CR}M\x T^{CR}M\to
TM/T^{CR}M$, the {\em real Levi--Form\/}. $(M,T^{CR}M)$ is called a
{\em partially integrable almost CR--manifold\/} if and only if $\Cal
L$ is non--degenerate and totally real, i.e.\ {}$\Cal L(J(\xi),J(\eta))=\Cal
L(\xi,\eta)$ for all $\xi,\eta\in T^{CR}M$, where $J$ denotes the
almost complex structure on $T^{CR}M$. In that case, choosing a local
trivialization of $TM/T^{CR}M$, $\Cal L$ is the imaginary part of a
Hermitian form. Here we consider the case where $n=2$, so $M$ has
dimension 5 and this Hermitian form is positive definite (for an
appropriate choice of the local trivialization). 

Typical examples of such manifolds are smooth hypersurfaces in a
six--dimensional smooth manifold $N$ endowed with an almost complex
structure $\tilde J$, which satisfy a non--degeneracy and an
integrability condition. In this case, we put $T^{CR}M=TM\cap\tilde
J(TM)$ and $J=\tilde J|_{T^{CR}M}$. To understand the non--degeneracy
and integrability conditions, it is more convenient to pass to
complexified tangent bundles. Since $T^{CR}M$ is a complex bundle, its
complexification $T^{CR}_\Bbb CM$ splits into a direct sum
$T_{1,0}M\oplus T_{0,1}M$ of a holomorphic and an antiholomorphic
part. Moreover, mapping $\xi,\eta\in\Ga(T_{1,0}M)$ to the class of
$-i[\xi,\bar\eta]$ defines a bundle valued Hermitian form $\Cal
L:T_{1,0}M\x T_{1,0}M\to T_\Bbb CM/T^{CR}_\Bbb CM=:QM$, the {\em Levi 
form\/}. The partial integrability condition from above is then
equivalent to the fact that $[\xi,\eta]\in\Ga(T^{CR}_\Bbb CM)$ for all
sections $\xi,\eta$ of $T_{1,0}M$, and the conditions of positive
definiteness is equivalent to $\Cal L$ being positive definite in
an appropriate local trivialization of $QM$. (Certainly, these
conditions also make sense for abstract almost CR manifolds). A
partially integrable almost CR manifold is called {\em integrable\/}
or a {\em CR--manifold\/} if and only if the subbundle $T_{1,0}M$ is
involutive. In particular, this is the case for hypersurfaces in
complex manifolds. 

By \cite[4.14]{CS}, 5--dimensional partially integrable almost
CR--manifolds are exactly the normal parabolic geometries
corresponding to $G=PSU(3,1)$ and the parabolic subalgebra of $\frak
g=\frak s\frak u(3,1)$ corresponding to the Dynkin diagram 
$\begin{picture}(50,8)
\put(5,2.5){\line(1,0){40}}
\put(5,2.5){\makebox(0,0){$\times$}}
\put(25,2.5){\makebox(0,0){$\bullet$}}
\put(45,2.5){\makebox(0,0){$\times$}}
\end{picture}%
$. Let us also consider $\tilde G=SU(3,1)$ and let $P$,
$G_0\subset G$,  or $\tilde P$, $\tilde G_0\subset \tilde G$ be the
corresponding subgroups as in \ref{2.3}. Then the semisimple part of
$\tilde G_0$ is  $SU(2)$ and the  center of $G_0$ is $\Bbb C$. 

In the Dynkin diagram notation, each
(complex) irreducible $\tilde G$-module $\Bbb V$ 
is given by the choice of three non--negative integers $a,b,c$
$$
\Bbb V=\dynkinG abc
.$$
More explicitly, $\dynkinG abc$ is the highest weight component in
$S^a{\Bbb C}^4\otimes S^b(\La^2{\Bbb C}^4)\otimes S^c({\Bbb C}^{4*})$,
where $S^i$ denotes the $i$--th symmetric power, 
and so these representations integrate to representation of $G$ if and
only if $a-c$ is congruent to $2b$ modulo four (the center of $\tilde
G$ consist of $\pm 1$ and $\pm i$ times the identity). 

The irreducible $\tilde P$--modules correspond to choices with $b$
non--negative while $a$ and $c$ may be arbitrary integers. Now, $b$
determines the representation of $SU(2)$ while the other two parameters
describe the action of the center of $\tilde G_0$.  We adopt the convention
used in \cite{BE}, i.e. the parameters give the linear combination of the
fundamental weights of $\tilde {\frak g}$ which is the highest weight of the
dual module to $\Bbb V$. In
this way, the resulting weights for our modules happen to be the same as
those in the dual pictures known from representation theory. For our
purposes, however, this has no importance and it is enough to say that the
distinguished two subbundles $T_{1,0}M$ and $ T_{0,1}M$ in the
complexified tangent space and the complexified quotient
$QM=T_\Bbb CM/T^{CR}_\Bbb CM$ have duals $T^*_{1,0}M$, $T^*_{0,1}M$
(quotients of the complexified cotangent bundle), and $Q^*M$, which
correspond to the modules 
$$
\Bbb T^*_{1,0}=\dynkin{-2}10, \quad \Bbb T^*_{0,1}=\dynkin01{-2}, \quad 
\Bbb Q^*=\dynkin{-1}0{-1}
.$$
Now, all $\tilde P$--modules are tensor products of symmetric powers 
$S^b(\Bbb T^*_{1,0})$ and suitable one-dimensional representations $\Bbb
E[a,c]$ corresponding to the Dynkin diagram $\dynkin a0c$. We shall write
$S^b(\Bbb T^*_{1,0})[a,c]$ for these modules and use the shorthand 
$\SSS bac$ for the corresponding natural bundles. In particular, 
\begin{gather*}
\SSS bac = S^b(\Bbb T^*_{1,0})[a,c]=\dynkin {a-2b}bc
\\
T^*_{0,1}=T^*_{1,0}[2,-2]=\SSS {1}2{-2}
\\
\SSS0{-1}{-1}=E[-1,-1]=Q^*
\\
\SSS0{-4}{0}=\La^2T^*_{1,0}\otimes Q^*.
\end{gather*}
Another important bundle is the dual to the kernel of the bilinear Levi form
$(\operatorname{ker}\Cal L)^*\subset T_{1,0}^*\otimes
T_{0,1}^*$ which corresponds to $\SSS22{-2}$.

All natural bundles $\SSS bac$ exist on manifolds $M$ with the so called
\emph{$SU(3,1)$--structures}, i.e. we have to choose coverings of the Cartan
$P$--bundle $\Cal G$ to the structure group $\tilde P$. This is clearly
equivalent to the choice of a fixed line bundle $E[1,0]$ such that its
fourth tensor power is $\La^2T_{1,0}M\otimes QM$. This is an analogous
situation to natural bundles and natural operators in conformal Riemannian
geometry which often depend on the choice of a spin structure.

\begin{figure}
$$
\xymatrix@C=-5mm{%
&&{\SSS b{a+2b}c}
\ar[dl]_{c+1}
\ar[dr]^{a+1}
\\
&{\SSS{b+c+1}{a+2b+2c+2}{-c-2}}
\ar[dl]_{b+1}
\ar[dr]^{a+1}
&&{\SSS{a+b+1}{a+2b}c}
\ar[dl]_{c+1}
\ar[dr]^{b+1}
\\
{\SSS c{a+b+2c+1}{-b-c-3}}
\ar[d]_{2a+2}
\ar[drr]|{\hbox to10mm{\hss\strut $\ssize a+b+2$\hss}}
&&{\SSS{a+b+c+2}{a+2b+2c+2}{-c-2}}
\ar[dll]|{\hbox to10mm{\strut\hss} }
\ar[d]^{2b+2}
\ar[drr]|{\hbox to10mm{\hss\strut $\ssize b+c+2$\hss}}
&&{\SSS a{a-b-3}{b+c+1}}
\ar[dll]|{\hbox to10mm{\strut\hss}}
\ar[d]^{2c+2}
\\
{\SSS c{b+2c}{-a-b-c-4}}
\ar[dr]_{b+1}
&&{\SSS {a+b+c+2}{a+b+2c+1}{-b-c-3}}
\ar[dl]^{a+1}
\ar[dr]_{c+1}
&&{\SSS a{a-b-c-4}{b}}
\ar[dl]^{b+1}
\\
&{\SSS{b+c+1}{b+2c}{-a-b-c-4}}
\ar[dr]_{c+1}
&&{\SSS{a+b+1}{a+b-c-2}{-b-2}}
\ar[dl]^{a+1}
\\
&&{\SSS b{b-c-3}{-a-b-3}}
}
$$
\caption{Bernstein--Gelfand--Gelfand sequences on partially integrable
5--dimensional almost CR structures}
\label{figure}
\end{figure}

Using the explicit description of the cohomology from 
Kostant's Bott--Borel--Weil theorem we obtain explicitly all 
natural bundles appearing in our main theorems. The computations are
done fairly simply in terms of the Dynkin diagram notation, see
\cite{BE} for the details. Furthermore, using elementary finite dimensional
representation theory one easily shows that there are no
homomorphisms between the semi--holonomic jet modules corresponding to
the items in the neighboring columns of the BGG--sequences, except
those which are indicated in Figure 
\ref{figure}. Let us also notice that the orders of the operators are easily
read off the homogeneities of the bundles with respect to the action of the
grading element in $G_0$ and the homogeneity of $\SSS bac$ is $a+c-b$. 
Thus we can summarize: 

\begin{thm} For each $SU(3,1)$--module $\Bbb V=\dynkinG abc$, the
BGG--sequence of invariant differential operators 
shown on Figure \ref{figure} exists on all 5--dimensional 
partially integrable almost CR--manifolds $M$  with a chosen
$SU(3,1)$--structure. The orders of the operators are indicated by the
labels over the arrows. Moreover, the sequence exists on all partially
integrable CR--manifolds if $a-2b+c \equiv 0\, \operatorname{mod}4$, 
and then all bundles in question can be constructed from $T^*_{1,0}M$
and $Q^*M$. If $M$ is flat, then the BGG--sequence is a complex which
computes the twisted de~Rham cohomology of $M$ with coefficients in
the bundle $VM$ corresponding to $\Bbb V$. 
\end{thm}

As a corollary, we immediately obtain

\begin{thm}
For all (integrable) 5--dimensional CR--manifolds, there is the resolution
of the sheaf of constant complex functions
$$
\xymatrix@R=2mm@C=4mm{%
&&&{\La^2T^*_{1,0}}
\ar[r]
\ar[ddr]
&{Q^*\otimes\La^2T^*_{1,0}}
\ar[dr]
\\
&&{T^*_{1,0}}
\ar[ur]
\ar[dr]
&&&{\otimes^2Q^*\otimes T^*_{1,0}}
\ar[dr]
\\
\Bbb C
\ar[r]
&{E[0,0]}
\ar[ur]
\ar[dr]
&& {(\operatorname{ker}\Cal L)^*}
\ar[uur]
\ar[r]
\ar[ddr]
&{Q^*\otimes(\operatorname{ker}\Cal L)^*}
\ar[ur]
\ar[dr]
&&\otimes^3Q^*
\\
&&{T^*_{0,1}}
\ar[ur]
\ar[dr]
&&&{\otimes^2Q^*\otimes T^*_{0,1}}
\ar[ur]
\\
&&&{\La^2T^*_{0,1}}
\ar[r]
\ar[uur]
&{Q^*\otimes\La^2T^*_{0,1}}
\ar[ur]
}
$$
which computes the de~Rham cohomology with complex coefficients. The
orders of the operators in the column in the middle of the diagram are
two,  while all the other ones are of first order.
\end{thm} 


This complex is a special instance of the so called Rumin complex on contact
geometries, \cite{Rum}, see also \cite{GarLee} for a refined version for the
CR--structures. In the homogeneous case, this complex was also mentioned in
\cite{BE}. Similar questions were also studied by Lychagin earlier, see e.g.
\cite{Lychagin} and the references therein. Notice that the dimensions of
the individual columns are 1, 4, 5, 5, 4, 1 (opposed to dimensions 1, 5, 10,
10, 5, 1 in the de~Rham complex).

\renewcommand{\thesection}{Appendix \Alph{section}}
\setcounter{section}{0}

\section{Infinite jets and Verma modules}\label{A}
\renewcommand{\thesection}{\Alph{section}}
The aim of this appendix is to provide differential geometers with basic
information on the links between jets and Verma modules, and in
particular to prove the correspondence between invariant differential
operators and homomorphisms of generalized Verma modules used in
\ref{2.6}.

\subsection{}\label{2.6a}
We have seen in \ref{2.6} that invariant operators 
$\Ga(E)\to\Ga(F)$ between homogeneous vector bundles over $G/P$ are
in bijective correspondence with $P$--homomorphisms $J^\infty(E)_o\to
F_o$, which factorize over some $J^r(E)_o$. 

First note that sections of $E$ can be identified with smooth
functions $G\to\Bbb E$, which are $P$--equivariant. Since this
identification is purely algebraic, it gives an identification of
infinite jets at $o$ of sections of $E$ with $P$--equivariant infinite
jets of smooth functions $G\to\Bbb E$ at $e\in G$. Now it is easy to
verify that in the picture of smooth equivariant functions, the action
of $G$ is given by $(g\cdot s)(g')=s(g^{-1}g')$. The corresponding
infinitesimal action of $\frak g$ is given by $(X\cdot s)(g)=-(R_X\cdot
s)(g)$, where $R_X$ denotes the {\em right--invariant\/} vector field
on $G$ generated by $X\in\frak g=T_eG$. For $X\in\frak p$, the
infinitesimal version of equivariancy of $s$ implies that $(X\cdot
s)(g)=X\cdot(s(g))$, but for general $X$ the value $(X\cdot s)(g)$
depends on the one--jet of $s$ at $g$. Thus we do not get an induced
action of $\frak g$ on finite jets, but for infinite jets we get a
well defined action of $\frak g$. Since this action is clearly
compatible with the action of $P$, it makes 
$J^\infty(E)_o$ into a $(\frak g,P)$--module. 

On the other hand, mapping each $X\in\frak g$ to the left invariant
vector field $L_X$ generated by $X$ induces an isomorphism between the
universal enveloping algebra $\Cal U(\frak g)$ and the algebra of
left invariant differential operators on $G$. Now we get a bilinear
map $J^\infty(E)_o\x(\Cal U(\frak g)\otimes\Bbb E^*)\to\Bbb K$ by
mapping $(j^\infty s(e),D\otimes\la)$ to $\la(D(s)(e))$, where $D$ is
a left invariant differential operator and $\la$ is an element of the
dual representation $\Bbb E^*$ to $\Bbb E$, and as above we view $s$
as an equivariant function on $G$. By equivariancy of $s$ this factors
to a bilinear map $J^\infty(E)_o\x(\Cal U(\frak g)\otimes_{\Cal
U(\frak p)}\Bbb E^*)\to\Bbb K$ because elements of $\Cal U(\frak p)$
act algebraically and this can be expressed as an action on $\la$. 

We claim that the above pairing is compatible with the actions of both
$\frak g$ and $P$. For the action of $\frak g$, let us take a typical
element $X_1\otimes\dots\otimes X_n\otimes\la\in \Cal U(\frak
g)\otimes_{\Cal U(\frak p)}\Bbb E^*$ and $X\in\frak g$. From above, we
know that $X\cdot j^\infty s(e)=-j^\infty (R_X\cdot s)(e)$. Pairing
this with $X_1\otimes\dots\otimes X_n\otimes\la$, we get 
$-\la((L_{X_1}\dots L_{X_n}\cdot R_X\cdot s)(e))$. Since left
invariant vector fields commute with right invariant ones, this equals
$-\la((R_X\cdot L_{X_1}\dots L_{X_n}\cdot s)(e))$. But this depends
only on $R_X(e)$, so we may as well replace $R_X$ by $L_X$, so this
coincides with $X\otimes X_1\otimes\dots\otimes X_n\otimes\la$
evaluated on $j^\infty s(e)$. 

The action of $b\in P$ on $\Cal U(\frak g)\otimes_{\Cal U(\frak p)}\Bbb
E^*$ is induced by mapping $D\otimes\la$ to $b\cdot D\otimes b\cdot
\la$, where $(b\cdot D)(s)=D(s\o r^{b^{-1}})\o r^b$ and $r^b$ denotes
the right multiplication by $b$. This obviously maps the anihilator of
the space of $P$--equivariant functions to itself and thus descends to
an action on $\Cal U(\frak g)\otimes_{\Cal U(\frak p)}\Bbb
E^*$. If $s$ is equivariant, then $(s\o r^{b^{-1}})(g)=b\cdot (s(g))$,
and thus $(b\cdot D)(s)(g)=b\cdot (D(s)(gb))$. But this implies that
pairing $j^\infty s(e)$ with $b\cdot D\otimes b\cdot\la$ we get
$(b\cdot\la)((b\cdot D)(s)(e))=\la (D(s)(b))$. On the other hand, the
action of $b$ on $J^\infty(E)_o$ is given by $b\cdot j^\infty
s(e)=j^\infty (s\o\ell_{b^{-1}})(e)$, where $\ell_b$ denotes the left
multiplication by $b$. Thus pairing $b^{-1}\cdot j^\infty s(e)$ with
$D\otimes\la$ we get $\la(D(s\o\ell_b)(e))$, which by left invariance
of $D$ coincides with $\la(D(s)(b))$. 

Now for any $k\in\Bbb N$, we have the natural projection
$J^\infty(E)_o\to J^k(E)_o$. On the other hand, the universal
enveloping algebra $\Cal U(\frak g)$ has a natural (infinite)
filtration $\Bbb K=\Cal U^0(\frak g)\subset\Cal U^1(\frak
g)\subset\dots$ such that $\Cal U(\frak g)=\cup_{i\in\Bbb N}\Cal
U^i(\frak g)$. In the picture of left invariant differential operators
on $G$, this is just the filtration by the order of operators. This
filtration clearly induces a filtration $\Cal F^i$ on $\Cal U(\frak
g)\otimes_{\Cal U(\frak p)}\Bbb E^*$, and each filtration component is
a $P$--submodule (but not a $\frak g$--submodule). The pairing of an
element of $\Cal F^k$ with an element $j^\infty(s)(e)\in
J^\infty(E)_o$ clearly depends only on $j^ks(e)$, so we get an induced
paring between $\Cal F^k$ and $J^k(E)_o$, and this induced pairing is
obviously non--degenerate and still compatible with the $P$--actions,
so since both sides are finite dimensional, they are dual
$P$--modules. 

Let us remark at this point that it is also possible to put
locally convex topologies on the spaces in question, such that they
become topologically dual $(\frak g,P)$--modules. Namely, one has to
view $J^\infty(E)_o$ as the limit of the system $\dots\to J^k(E)_o\to
J^{k-1}(E)_o\to\dots$, while $\Cal U(\frak g)\otimes_{\Cal U(\frak
p)}\Bbb E^*$ has to be topologized as a direct sum of finite
dimensional spaces. 

\subsection{}\label{2.7}
Let $\Bbb E$ and $\Bbb F$ be $P$--representations, $E$ and $F$ the
corresponding bundles and $\ph:J^k(E)_o\to F_o=\Bbb F$ a
$P$--homomorphism. By the duality shown above, we can view the dual
map $\ph^*$ as a $P$--homomorphism $\Bbb F^*\to\Cal F^k\subset \Cal
U(\frak g)\otimes_{\Cal U(\frak p)}\Bbb E^*$. Conversely, if we have a
$P$--homomorphism $\Bbb F^*\to \Cal U(\frak
g)\otimes_{\Cal U(\frak p)}\Bbb E^*$, then this has values in some
$\Cal F^i$ since $\Bbb F^*$ is finite dimensional, so dualizing it
corresponds to a $P$--homomorphism $J^i(E)_o\to F_o$. Consequently,
we see that the space of invariant operators $\Ga(E)\to\Ga(F)$ is
isomorphic to $\Hom_P(\Bbb F^*,\Cal U(\frak g)\otimes_{\Cal
U(\frak p)}\Bbb E^*)$.

By Frobenius reciprocity 
the latter space is isomorphic to  
$$\Hom_{(\frak g,P)}(\Cal U(\frak
g)\otimes_{\Cal U(\frak p)}\Bbb F^*,\Cal U(\frak g)\otimes_{\Cal
U(\frak p)}\Bbb E^*).$$ 
This isomorphism is quite simple to prove: If $\ph:\Bbb F^*\to\Cal U(\frak
g)\otimes_{\Cal U(\frak p)}\Bbb E^*$ is a $P$--homomorphism, then
$\tilde\Ph(A\otimes \la)=A\cdot\ph(\la)$ defines a $(\frak
g,P)$--homomorphism $\Cal U(\frak
g)\otimes\Bbb F^*\to\Cal U(\frak g)\otimes_{\Cal U(\frak p)}\Bbb E^*$,
and since $\ph$ is a $P$--homomorphism, this factors to a $(\frak
g,P)$--homomorphism $\Ph$ between the required spaces. Conversely, we
put $\ph(\la)=\Ph(1\otimes\la)$ and this clearly is a
$P$--homomorphism if $\Ph$ is a $(\frak g,P)$--homomorphism. 

\renewcommand{\thesection}{Appendix \Alph{section}}
\section{Adjointness of $\partial$ and $\partial^*$}\label{B}
\renewcommand{\thesection}{\Alph{section}}
\subsection{}\label{4.2}
As promised in the beginning of Section \ref{4}, we show that the 
operators $\partial$ and
$\partial^*$ are adjoint operators with respect to a certain inner
product on $C^n(\frak g_-,\Bbb V)$. To construct this inner product,
we have to distinguish between the real and the complex case. Let us
start with the case where $\frak g$ and $\Bbb V$ are complex. Since
the grading element $E\in\frak g_0$ is semisimple, we can find a
Cartan subalgebra $\frak h\subset\frak g$ which contains $E$. Then
each root space for this Cartan subalgebra is contained in some $\frak
g_i$. Let $\frak u$ be a compact real form of $\frak g$ with a Cartan
subalgebra $\frak h_0$ contained in $\frak h$, and let $\si$ be the
complex conjugation with respect to this real form. By definition of
$E$, the map $\ad(E)\o\ad(E)$ acts on $\frak g_i$ by multiplication by
$i^2$, so for the Killing form we have $B(E,E)>0$. Consequently, we
must have $\si(E)=-E$, and thus $\si(\frak g_i)=\frak g_{-i}$ for all
$i=-k,\dots, k$. Now one immediately verifies directly that
$B^*(X,Y):=-B(X,\si(Y))$ is a positive definite Hermitian inner product
on $\frak g$, such that the decomposition $\frak g=\frak
g_{-k}\oplus\dots\oplus\frak g_k$ is an orthogonal direct sum. In
particular, this induces a Hermitian inner product on $\frak g_-$. 

Next, since $\frak u$ is a compact real form, there is a positive
definite Hermitian inner product $\langle\ ,\ \rangle$ on $\Bbb V$ such
that the elements of $\frak u$ act as skew--Hermitian operators. But
this immediately implies that for each $X\in\frak g$ and
$v_1,v_2\in\Bbb V$, we have $\langle X\cdot v_1,v_2\rangle=-\langle
v_1,\si(X)\cdot v_2\rangle$. Together with the inner product on $\frak
g_-$ constructed above we get a positive definite Hermitian inner
product on $C^n(\frak g_-,\Bbb V)$ for each $n$.

In the real case, the situation is slightly more complicated. In this
case we have to construct appropriate involutions $\si$ on the
individual simple factors separately, and we have to distinguish
between the case where the complexification of a simple factor is
again simple and the case where it is not. Note that the simple
factors of a $|k|$--graded Lie algebra are themselves $|\ell|$--graded
for some $\ell\leq k$ and that the grading element of $\frak g$ is
just the sum of the grading elements of the simple factors. 

If we have a real simple algebra $\frak g$ whose complexification is
not simple, then it is well known that $\frak g$ is actually the
underlying real Lie algebra of a complex simple Lie algebra. In this
case, we can proceed as above to get a compact real form $\frak
u\subset\frak g$ and the corresponding involution $\si$.  

In the case where both $\frak g$ and its complexification $\frak
g^\Bbb C$ are simple, we choose a Cartan subalgebra $\frak
h\subset\frak g^\Bbb C$ which contains the element $E\in\frak g$. By
\cite[Expos\`e 11, Th\'eor\`eme 3]{SemSophLie} there is a compact real
form $\frak u\subset\frak g^\Bbb C$ with Cartan subalgebra $\frak
h_0\subset\frak h$ such that the complex conjugation $\si$ with
respect to $\frak u$ commutes with the complex conjugation with
respect to $\frak g$, and thus $\si(\frak g)=\frak g$.

The involutions on the simple factor together define an involution of
$\frak g$ and as above one uses the Killing form on $\frak g$ and
$\si$ to get a positive definite inner product on $\frak g$ and on
$\frak g_-$. If the representation $\Bbb V$ is not already complex,
then we can pass to its complexification to get a Hermitian inner
product such that $\langle X\cdot v_1,v_2\rangle=-\langle
v_1,\si(X)\cdot v_2\rangle$ as above, an in both cases the real part
of this Hermitian product gives a positive definite inner product on
$\Bbb V$ which we use together with the inner product on $\frak g_-$
to get a positive definite inner product on $C^n(\frak g_-,\Bbb V)$. 

\begin{prop}\label{4.3}
The differential $\partial:C^n(\frak g_{-1},\Bbb V)\to C^{n+1}(\frak
g_{-1},\Bbb V)$ and the codifferential 
$\partial^*: C^{n+1}(\frak g_{-1},\Bbb V)\to C^n(\frak g_{-1},\Bbb V)$ are
adjoint operators with respect to the inner products constructed in
\ref{4.2} above.
\end{prop}
\begin{proof}
The point about this is that in each case the inner product of
$f_1,f_2\in C^n(\frak g_-,\Bbb V)$ can be computed as 
$\Cal F(f_2)(f_1)$, where $\Cal F$ is a linear (over the reals)
isomorphism $C^n(\frak g_-,\Bbb V)\to C^n(\frak p_+,\Bbb V^*)$. The map
$\Cal F$ is defined by $\Cal F(f)(Z_1,\dots,Z_n)(v):=\langle
f(\si(Z_1),\dots,\si(Z_n)),v\rangle$ for $Z_i\in\frak p_+$ and $v\in
\Bbb V$, where $\si$ is the involution constructed in \ref{4.2} and the
inner product is in $\Bbb V$. But then the compatibility of the inner
product on $\Bbb V$ with the action of $\frak g$ implies that $\Cal
F(\partial(f))=\partial(\Cal F(f))$. Thus we can compute:
\begin{multline*}
\langle\partial^*(f_1),f_2\rangle=\Cal F(f_2)(\partial^*(f_1))=
\partial(\Cal F(f_2))(f_1)\\
=\Cal F(\partial(f_2))(f_1)=
\langle f_1,\partial(f_2)\rangle
\end{multline*}
\end{proof}

\end{document}